\documentclass[12pt]{amsart}
\usepackage{amsfonts}
\usepackage{amssymb}
\usepackage{color}
\usepackage{hyperref}
\numberwithin{equation}{section}
\theoremstyle{plain}
 \newtheorem{thm}{Theorem}[section]
 \newtheorem{lem}[thm]{Lemma}
 
 \newtheorem{prop}[thm]{Proposition}
 
\theoremstyle{definition}
\newtheorem{defn}[thm]{Definition}
 
 \newtheorem{rem}[thm]{Remark}

\newcommand{\al}{\alpha}
\newcommand{\bt}{\beta}
\newcommand{\gm}{\gamma}

\newcommand{\dl}{\delta}

\newcommand{\ep}{\varepsilon}

\newcommand{\ld}{\lambda}

\newcommand{\q}{\quad}

\newcommand{\wh}{\widehat}
\newcommand{\wt}{\widetilde}
\newcommand{\la}{\langle}
\newcommand{\ra}{\rangle}

\newcommand{\R}{\mathbb{R}}

\newcommand{\rd}{\mathbb R ^d}
\newcommand{\law}{\mathcal L}
\newcommand{\sek}{\int_0^{\infty}}
\newcommand{\n}{\noindent}

\makeatletter
\renewcommand\subsection{\@startsection{subsection}{2}%
  \z@{.5\linespacing\@plus.7\linespacing}{.5\linespacing}%
  {\normalfont\scshape\centering}}
\makeatother

\begin{document}
\setlength{\baselineskip}{18pt}
\setlength{\parindent}{1.8pc}

\vskip 5mm

\title{Classes of infinitely divisible distributions on $\rd$ related to the class of
selfdecomposable distributions}\maketitle

\begin{center}
{Makoto Maejima\footnote{Corresponding author: maejima@math.keio.ac.jp}, Muneya Matsui and Mayo Suzuki\\
\vskip 3mm
Department of Mathematics, Keio University\\
3-14-1, Hiyoshi, Kohoku-ku\\
Yokohama 223-8522, Japan}
\end{center}

\vskip 10mm
\begin{center}
ABSTRACT
\end{center}
This paper studies new classes of infinitely divisible distributions on $\rd$.
Firstly, the connecting classes with a continuous parameter
between the Jurek class and the class of selfdecomposable distributions are revisited.
Secondly, the range of the parameter is extended to construct new classes and characterization
in terms of stochastic integrals with respect to L\'evy processes are given.
Finally,  the nested subclasses of those classes are discussed and characterized in two ways:
One is by stochastic integral representations and another is in terms of L\'evy measures.

\vskip 10mm


\allowdisplaybreaks
\section{Introduction}

Let $I(\rd)$ be the class of all infinitely divisible distributions on $\rd$ and
$I_{\log}(\rd)=\{\mu\in I(\rd) : \int_{|x|>1}\log|x|\mu(dx)<\infty\}$.
Let $\wh \mu (z), z\in \rd$, be the characteristic function of $\mu\in I(\rd)$.

In this paper, we first revisit the classes in $I(\rd)$ connecting the class of selfdecomposable distributions 
($L(\rd)$, say) and the Jurek class (the class of $s$-selfdecomposable distributions, 
($U(\rd)$, say), see Jurek (1985)). 
Those connecting classes were already studied by O'Connor (1979) in $I(\R ^1)$ and by Jurek (1988) in $I(E)$, 
where $E$ is a Banach space.
Throughout this paper, we treat the case $I(\rd)$.
Although there are several equivalent definitions of $L(\rd)$ and $U(\rd)$, we use here their definitions in 
terms of L\'evy measures.
Then we study more general classes including the classes above and nested subclasses of those classes.

The L\'evy-Khintchine representation of $\wh\mu$ we use in this paper is
$$
\wh\mu(z) = \exp \biggl\{ -{2}^{-1}\la z,Az \ra + i\la \gm,z \ra +\int_{\rd}\left(e^{i\la z,x \ra}-1
-\frac{i\la z,x \ra}{1+|x|^2}\right)\nu(dx) \biggr\},
$$
where $A$ is a nonnegative-definite symmetric $d\times d$ matrix, $\gm\in\rd$ and $\nu$ is the L\'evy measure
satisfying $\nu(\{0\})=0$ and $\int_{\rd}(|x|^2\wedge 1)\nu(dx)<\infty$.
We call $(A,\nu, \gm)$ the L\'evy-Khintchine triplet of $\mu$ and we write $\mu=\mu_{(A,\nu,\gm)}$
when we want to emphasize the triplet.

The polar decomposition of the L\'evy measure $\nu$ of $\mu\in I(\rd)$, with $0<\nu(\R^d)\le\infty$, 
is the following:
There exist a measure $\ld$ on $S=\{\xi\in\rd : |\xi |=1\}$
with $0<\ld(S)\le\infty$ and
a family $\{\nu_{\xi}\colon \xi\in S\}$ of measures on $(0,\infty)$ such that
$\nu_{\xi}(B)$ is measurable in $\xi$ for each $B\in\mathcal B((0,\infty))$,
$0<\nu_{\xi}((0,\infty))\le\infty$ for each $\xi\in S$ and
\begin{align}
\nu(B)=\int_S \ld(d\xi)\int_0^{\infty} 1_B(r\xi)\nu_{\xi}(dr), \,\,
B\in \mathcal B (\mathbb R^d \setminus \{ 0\}).
\end{align}
Here $\ld$ and $\{\nu_{\xi}\}$ are uniquely determined by $\nu$
up to multiplication of a measurable function $c(\xi)$ and $\frac1{c(\xi)}$, respectively,
with $0<c(\xi)<\infty$.
We say that $\mu$ or $\nu$ has the polar decomposition $(\ld ,\nu_{\xi})$ and
$\nu_{\xi}$ is called the radial component of $\nu$.
(See, e.g., Barndorff-Nielsen et al. (2006), Lemma 2.1.)

The connecting classes between $U(\rd)$ and $L(\rd)$ mentioned above are also 
characterized by mappings with a parameter from $I(\rd)$  into $I(\rd)$.
We extend the range of the parameter and first study the classes defined by these mappings.
These mappings are the special cases studied in Sato (2006b) as will be mentioned later.

We start with following classes, where the classes $U(\rd)$ 
and $L(\rd)$ are two known special classes.

\begin{defn}
(Class $K_{\al}(\rd)$)  
Let $\al < 2$. 
We say that $\mu\in I(\rd)$ belongs to the class $K_{\al}(\rd)$ if $\nu =0$ or $\nu \ne 0$ and,
in case $\nu\ne 0$, $\nu_{\xi}$ in (1.1) satisfies
\begin{equation}
\nu_{\xi}(dr) = r^{-\al -1}\ell_{\xi}(r)dr , \,\, r>0,
\end{equation}
where $\ell_{\xi}(r)$ is nonincreasing in $r\in (0,\infty)$ for $\ld$-a.e.\,$\xi$ and is measurable in $\xi$
for each $r>0$, and $\lim_{r\to\infty}\ell_{\xi}(r) = 0$.
\end{defn}

\begin{rem}
(i) Because of the condition that $\lim_{r\to\infty}\ell_{\xi}(r)=0$, $K_{\al}(\rd), 0<\al<2$, does not include
the class of $\al$-stable distributions, but does include the class of $(\al +\ep)$-stable distributions for
any $\ep\in (0,2-\al)$.
It also includes tempering $\al$-stable distributions,
which are defined by (1.2) with a completely monotone function $\ell_{\xi}(r)$ on $(0,\infty)$
such that $\lim_{r\to 0}\ell_{\xi}(r) = 1$ and $\lim_{r\to\infty}\ell_{\xi}(r) = 0$.
(See Rosi\'nski (2007).) \\
(ii) Let $\nu$ be the L\'evy measure of  $\mu\in I(\rd)$ and $\al >0$.
Since $\int_{\rd}|x|^{\dl}\mu(ds)<\infty$ if and only if $\int_{|x|>1}|x|^{\dl}\nu(dx)<\infty$,
(see, e.g. Sato (1999) Theorem 25.3,) $\mu\in K_{\al}(\rd)$ has the finite $\dl$-moment
for any $0<\dl<\al$.
This fact is the same as for $\al$-stable distributions.
\end{rem}

\begin{rem}
(i) The Jurek class 
$ U(\rd)$ is $K_{-1}(\rd)$ and 
the class of selfdecomposable distributions
$L(\rd)$ is $K_0(\rd)$. \\
(ii) Let $ \al <\bt < 2$.  Then
$
K_{\bt}(\rd) \subset K_{\al}(\rd).
$
This is trivial from the definition. 
\end{rem}
Therefore, $K_{\al}(\rd), -1\le \al\le 0$, are connecting classes with a continuous parameter $\al$
between the classes $U(\rd)$ and $L(\rd)$, as mentioned in the beginning of this section.

This paper is organized as follows.
In Section 2, some known results related to the classes $K_{\al}(\rd)$ are mentioned.
In Section 3, we give a complete proof for the decomposability of the distributions in $K_{\al}(\rd), \al <0$.
In Section 4, we define mappings $\Phi_{\al}, \al \in \R,$ in terms of stochastic integrals with respect to L\'evy processes
related to the classes $K_{\al}(\rd)$ and determine those domains and ranges.
The proofs for the ranges are given in Section 5.
In Section 6, we construct nested subclasses of the ranges of $\Phi_{\al}$ by iterating the mapping $\Phi_{\al}$.
Then we firstly determine the domains $\mathfrak D (\Phi_{\al}^{m+1}), m=1,2,3,...$, and secondly characterize 
the ranges of the mappings $\Phi_{\al}^{m+1}$ in two ways:
One is by stochastic integral representations and another is in terms of L\'evy measures.

\vskip 10mm
\section{Known results}

In this section, we explain several results from O'Connor (1979) and Jurek (1988).

\n
{\bf 1.}
(Characterization by the decomposability.)

O'Connor (1979) defined the classes $K_{\al}(\R ^1), -1<\al<0$, as in Definition 1.1, and proved that
$\mu\in K_{\al}(\R ^1)$ if and only if for any $c\in (0,1)$ there exists $\mu_c\in I(\R^1)$ such that
\begin{equation}
\label{eq:self-decomp}
\wh \mu (z) = \wh \mu (cz)^{c^{-\al}}\wh\mu _c(z).
\end{equation}
His proof used L\'evy measures.
However, his proof for getting the convexity of L\'evy density on $(-\infty,0)$ and the concavity on 
$(0,\infty)$ (in the proof of his Theorem 3 in O'Connor (1979)) is not clear to the authors of this
paper.
So, we will give our proof in Section 3, extending the range of $\al$ to $(-\infty, 0)$.

Jurek (1988) defined the classes $U_{\al}(E), -1\le\al\le 0$, where $E$ is a Banach space, as the classes of
limiting distributions as follows. 
$\mu\in  U_{\al}(E)$ if and only if there exists a sequence $\{ \mu_j\}\subset I(E)$ such that
\begin{equation}
\label{eq:self-decomp-limit}
\lim_{n\to\infty}n^{-1}(\mu_1*\mu_2*\cdots*\mu_n)^{*n^{\al}}=\mu,
\end{equation}
where $(n^{-1}(\mu_1*\mu_2*\dots*\mu_n)^{n^{-1}})(B):=(\mu_1*\mu_2*\dots*\mu_n)^{n^{-1}}(nB),B\in\mathcal B (E)$.
He then showed the decomposability (\ref{eq:self-decomp}) as a
consequence of (\ref{eq:self-decomp-limit}).
So, as a result, we see that $K_{\al}(\rd)= U_{\al}(\rd)$, but there is no proof by using L\'evy measures in Jurek (1988).
This is another reason why we will give our proof in Section 3.
Our proof will use L\'evy measures in the same way as in the proof of Theorem 15.10 of Sato (1999) for 
selfdecomposability.
\vskip 3mm
\n
{\bf 2.}
(Characterization by the stochastic integrals with respect to L\'evy processes.)

Let $-1\le \al < 0$.
Jurek (1988) showed that
$\mu\in U_{\al}(E)$ if and only if there exists a L\'evy process $\{X_t\}$ on $E$ such that
\begin{equation}
\label{eq:law-Jurek}
\mu =\law \left ( \int_0^1 t^{-1/\al}dX_t\right ),
\end{equation}
where $\law (X)$ is the law of a random variable $X$.
For the case $\al =0$, the following is known (Wolfe (1972) and others).
$\mu\in K_{0}(\rd)$ if and only if there exists a L\'evy process $\{X_t\}$ on $\rd$ 
satisfying $E[\log ^+ |X_1|] <\infty$ such that
$$
\mu =\law \left ( \sek e^{-t}dX_t\right ).
$$

\begin{rem} (\ref{eq:law-Jurek}) can have a different form.
Change the variables from $t$ to $s$ by $t=1+\al s$.
Then
$$
\mu = \law\left (-\int_0^{-1/\al}(1+\al s)^{-1/\al}dX_{1+\al s}\right).
$$
If we define another L\'evy process $\{ \wt X_t\}$ by $\wt X_s= -X_{1+\al s}$, then we have
\begin{equation}
\label{eq:law-alpha}
\mu = \law\left (\int_0^{-1/\al}(1+\al s)^{-1/\al}d\wt X_s\right).
\end{equation}
(\ref{eq:law-alpha}) will be seen in Definition 4.1 with $\al <0$ below and this expression is more natural when we 
consider the case $\al =0$ as
we will see in Remark 4.7 later.
\end{rem}

\vskip 10mm

\section{Decomposability of distributions in $K_{\al}(\rd)$}

As mentioned before, the classes $U(\rd)$ and $L(\rd)$ have characterizations in terms of characteristic functions.
Namely,
$\mu\in U(\rd)$ if and only if
for any $c\in (0,1)$, there exists $\mu_c(z)\in I(\rd)$ such that
$$
\wh{\mu}(z) = \wh{\mu} (cz)^c\wh{\mu}_c(z),
$$
and
$\mu\in L(\rd)$ if and only if
for any $c\in (0,1)$, there exists $\mu_c(z)\in I(\rd)$ such that
$$
\wh{\mu}(z) = \wh{\mu} (cz)\wh{\mu}_c(z).
$$

As we announced in Section 2, we give our proof of characterization of
 $K_{\al}(\rd)$ in a similar way as follows.

\begin{thm}
Let $\al <0$.
$\mu\in K_{\al}(\rd)$ if and only if
for any $c\in (0,1)$, there exists $\mu_c\in I(\rd)$ such that
$$
\wh{\mu}(z) = \wh{\mu} (cz)^{c^{-\al}}\wh{\mu}_c(z).
$$
\end{thm}

\vskip 3mm
\n
{\it Proof.}
(\lq\lq Only if '' part.)  
It is enough to consider the case with $A=O$ and $\gm =0$.
Suppose $\mu\in K_{\al}(\rd)$ and the polar decomposition of the 
L\'evy measure of $\mu$ is $(\ld, \nu_{\xi})$, with $\nu_{\xi}(dr) = r^{-\al-1}{\ell}_{\xi}(r)dr$.
Then we have
\begin{align*}
\wh\mu (z) 
& = \exp \left \{ \int _S \ld (d\xi) \sek \left ( e^{i\la z, r\xi\ra} - 1- \frac{i\la z,r\xi\ra}{1+r^2}
\right ) \frac1{r^{\al+1}}\ell_{\xi}(r)dr\right \}.
\end{align*}
Thus,
\begin{align*}
\wh\mu& (cz)^{c^{-\al}} \\
& = \exp \left \{ c^{-\al}\int_S \ld(d\xi) \sek
\left ( e^{i\la z, cr\xi\ra} - 1- \frac{i\la z,cr\xi\ra}{1+r^2}
\right )\frac{1}{r^{\al+1}}\ell_{\xi}(r)dr\right \}\\
& = \exp \left \{  \int _S \ld (d\xi) \sek \left ( e^{i\la z, u\xi\ra} - 1- 
\frac{i\la z,u\xi\ra}{1+(u/c)^2}
\right ) \frac{1}{u^{\al+1}}\ell_{\xi}\left (\frac{u}{c}\right)du\right \}\\
& = \wh{\mu}(z) 
\exp \left \{ \int _S \ld (d\xi) \sek{ i\la z,u\xi\ra} \left (\frac1{1+u^2} - \frac1{1+(u/c)^2}
\right ) \frac1{u^{\al+1}}\ell_{\xi}\left(\frac{u}{c}\right ) du\right \}\\
& \hskip 10mm\times 
\exp\left\{ -\int_S \ld(d\xi)\sek \left (e^{i\la z, u\xi\ra} - 1 - \frac{i\la z, u\xi\ra}{1+u^2}\right )\frac1{u^{\al+1}}
\left (\ell_{\xi} (u)- \ell_{\xi} \left (\frac{u}{c}\right) \right )du \right \}\\
& =:  \wh{\mu}(z)  e^{i\la z, a_c\ra}(\rho_c(z))^{-1},
\end{align*}
where
$$
a_c = \int _S \xi \ld (d\xi) \sek u^{-\al}\left (\frac1{1+u^2} - \frac1{1+(u/c)^2}
\right ) \ell_{\xi}\left(\frac{u}{c}\right )du
$$
and
$$
\wh\rho _c(z) = \exp\left\{ \int_S \ld(d\xi)\sek \left (e^{i\la z, u\xi\ra} - 1 - \frac{i\la z, u\xi\ra}{1+u^2}\right )\frac 1{u^{\al+1}}
\left (\ell_{\xi}(u) - \ell_{\xi }\left (\frac{u}{c}\right )\right)du \right \}.
$$
We have to check the finiteness of $a_c$ and that $\rho_c\in I(\rd)$.

Since $\nu$ is a L\'evy measure, we have
$\int_S\ld(d\xi)\int_0^\infty(r^2\wedge1)\nu_\xi(dr)<\infty$, which implies 
\begin{equation*}
\int_S\ld(d\xi)\int_0^1 r^{-\al+1}\ell_\xi(r)dr<\infty\ \text{and}\ \int_S\ld(d\xi)\int_1^\infty r^{-\al-1}\ell_\xi(r)dr<\infty.
\end{equation*}
Furthermore, this concludes
\begin{align*}
|a_c| &\le \int_S|\xi|\ld(d\xi) \int_0^\infty u ^{-\al }\left|\frac{1}{1+u^2}-\frac{1}{1+({u}/{c})^2}
\right|\ell_\xi({u}/{c})du\\
&= \int_S\ld(d\xi) \int_0^\infty c^{1-\al}v^{-\al}
 \left|\frac{1}{1+(cv)^2}-\frac{1}{1+v^2}\right|\ell_{\xi}(v)dv \\
&= \int_S\ld(d\xi) \int_0^\infty 
\frac{c^{1-\al} v^{1-\alpha}}{(1+(cv)^2)}  \ell_{\xi}(v)dv \\
& \le c^{1-\al}\int_S\ld(d\xi)\int_0^1v^{1-\al}\ell_\xi(v)dv  + c^{1-\al}\int_S\ld(d\xi)
\int_1^\infty \frac{v^{1-\al}}{1+(cv)^2}\ell_\xi(v)dv <\infty.
\end{align*}
This shows the finiteness of $a_c$.

With respect to $\rho_c$,
since $0<c<1$ and $\ell_{\xi}$ is nonincreasing, we have $h_\xi(u):= u^{-\al-1}(\ell_\xi(u) - \ell_{\xi}(u/c)) \geq 0.$  
Thus,  $\nu_{\rho}(B) = \int_S \ld (d\xi)\sek 1_B(r\xi)h_\xi(r)dr$ is a nonnegative measure.
Furthermore, we have
$$
\int_0^\infty(r^2\wedge1) \nu_\rho(dr)
= \int_S\ld(d\xi)\int_0^\infty(r^2\wedge1) r^{-\al-1}(\ell_\xi(r) - \ell_{\xi}(r/c))dr
< \infty ,
$$
because $\int_S\ld(d\xi)\int_0^\infty(r^2\wedge1)r^{-\al -1}\ell_{\xi}(r)dr<\infty.$
Therefore, $\nu_{\rho}$ is a L\'evy measure, and $\rho_c\in I(\rd) $ by the uniqueness of 
L\'evy-Khintchine representation.
Thus, if we put $\wh\mu _c(z) = \wh\rho_c (z) e^{-i\la z, a_c\ra}$, we have
$$
\wh{\mu}(z) = \wh{\mu} (cz)^{c^{-\al}}\wh{\mu}_c(z).
$$
\lq\lq Only if" part is now proved.

(\lq\lq If" part.)
Conversely, suppose that $\mu\in I(\rd)$ satisfies that
for any $c\in (0,1)$, there exists $\mu_c(z)\in I(\rd)$ such that
$\wh{\mu}(z) = \wh{\mu} (cz)^{c^{-\al}}\wh{\mu}_c(z).$
Since
$$\wh\mu(z) = \exp \biggl\{ -{2}^{-1}\la z,Az \ra + i\la \gm,z \ra +\int_{\rd}\left(e^{i\la z,x \ra}-1
-\frac{i\la z,x \ra}{1+|x|^2}\right)\nu(dx) \biggr\},$$
we have
\begin{align*}
\wh\mu(cz)^{c^{-\al}} 
&= \exp \biggl\{ -{2}^{-1}c^{-\al}\la cz,Acz \ra + ic^{-\al}\la \gm,cz \ra +\int_{\rd}
\left(e^{i\la cz,x \ra}-1-\frac{i\la cz,x \ra}{1+|x|^2}\right)c^{-\al}\nu(dx) \biggr\}\\
&= \exp \biggl\{ -{2}^{-1}\la z,c^{2-\al}Az \ra + i\la c^{1-\al}\gm,z \ra \\
& \hskip20mm +\int_{\rd}\left(e^{i\la z,y \ra}-1-\frac{i\la z,y \ra}{1+|y|^2}\right)c^{-\al}\nu 
\left(\frac{dy}{c}\right) \\
& \hskip20mm +i\left \la z,\int_{\rd}y\left(\frac{1}{1+|y|^2} - \frac{1}{1+|{y}/{c}|^2}\right)c^{-\al}\nu 
\left(\frac{dy}{c}\right) 
\right \ra\biggr\}.
\end{align*}
Since $\mu \in I(\rd),  \wh\mu(z)\ne 0$ for any $z\in\rd$.  Then we have
$$ \wh\mu_c(z)=\exp \biggl\{-2^{-1}\la z,A_c z \ra + i\la \gm_c ,z\ra+ \int_{\rd}\left(e^{i\la z,x \ra}-1-
\frac{i\la z,x \ra}{1+|x|^2}\right)\left(\nu(dx)-c^{-\al}\nu \left(\frac{dx}{c}\right)\right)  \biggr\},$$
where
$
A_c=(1-c^{2-\alpha})A
$
and
$$
\gm_c = (1- c^{1-\al})\gm  - \int_{\rd}y\left(\frac{1}{1+|y|^2} 
- \frac{1}{1+|{y}/{c}|^2}\right)c^{-\al}\nu \left(\frac{dy}{c}\right ).
$$
Since $\mu_c\in I(\rd)$,
$\nu^c(B):=\nu(B) - c^{-\al}\nu(c^{-1}B)$ is a L\'evy measure for any $c\in (0,1)$.
Recall that the polar decomposition of $\nu$ is $\nu(B) = \int_S\ld(d\xi)\int_0^\infty 1_B(r\xi)\nu_\xi(dr).$  
Then,
\begin{align*}
\nu^c(B)
&= \int_S\ld(d\xi)\int_0^\infty(1_B(r\xi)\nu_\xi(dr)-1_{c^{-1}B}(r\xi)c^{-\al}\nu_\xi(dr))\\
&= \int_S\ld(d\xi)\int_0^\infty 1_B(r\xi) \left(\nu_\xi(dr)-c^{-\al}\nu_\xi \left(\frac{dr}{c}\right) \right).
\end{align*}
It remains to show that 
$$
\nu_{\xi}(dr) = r^{-\al-1}\ell_{\xi}(r)dr
$$
for some nonincreasing function $\ell_{\xi}$ measurable in $\xi$ . 
For that, we consider a measure $r^\alpha\nu_\xi(dr)$ on $(0,\infty)$ and let 
$$
H_\xi (x):=\int_{e^{-x}}^\infty r^\alpha \nu_\xi(dr).
$$
Here $H_\xi (x)$ is measurable in $\xi$.
We also put 
\begin{align*}
H^c_\xi(x):&=H_\xi(x)-H_\xi(x+\log c) \\
&= \int_{e^{-x}}^\infty r^\alpha\nu_\xi(dr)-\int_{e^{-x}/c}^\infty r^\alpha\nu_\xi(dr) \\
&= \int_{e^{-x}}^{\infty} r^{\alpha}\left(
\nu_\xi(dr)-c^{-\alpha}\nu_\xi\left(\frac{dr}{c} \right)
\right). 
\end{align*}
Since $\nu^c(dr)$ is a L\'evy measure, $H_\xi^c(x)$ is nonnegative and is
nondecreasing for $\lambda-$almost every $\xi$. Moreover, $H_\xi(x)$ is convex on $(-\infty,\infty)$ as shown below.

Let $s\in\mathbb{R}, u>0$ and $c\in(0,1)$. Then $H_\xi^c(s+u)\ge H_\xi^c(s)$, and thus 
$$ 
H_\xi(s+u)-H_\xi(s+u+\log c)\ge H_\xi(s)-H_\xi(s+\log c)\ge 0, 
$$
which is 
\begin{equation}
H_\xi(s+u)-H_\xi (s) \ge H_\xi(s+u+\log c)-H_\xi(s+\log c)\ge 0. \label{ineq:continuous-1}
\end{equation}
Then $H_\xi$ is convex for $\lambda-$almost every $\xi$, as in Sato
(1999) pp. 95--96. Furthermore, repeating the argument in p.96 of Sato
(1999) we can write 
$$
H_\xi (x)=\int_{-\infty}^x h_\xi (t)dt,
$$
where $h_\xi(t)$ is some left-continuous nondecreasing function in
$u$. Hence $h_x(t)$ is measurable in $\xi$.
Now put
$$
 H_\xi (-\log x)= \int_{-\infty}^{-\log x}h_\xi(t)dt =\int_{x}^\infty h_\xi(-\log r)r^{-1}dr,
$$
then, the definition of $H$, we have
$$
\int_x^\infty r^\alpha \nu_\xi(dr)=\int_x^\infty h_\xi (-\log r)r^{-1}dr,
$$
which implies
$$
\nu_\xi(dr)=r^{-\alpha-1}h_\xi(-\log r)dr.
$$
Since $h_\xi$ is nondecreasing, we have $h_\xi(-\log r)$ is a nonincreasing function, and putting
$\ell_{\xi}(r)=h_\xi(-\log r)$, we complete the proof. 
\qed

\vskip 10mm

\section{Mappings defined by stochastic integrals related to $K_{\al}(\rd)$}

We are now going to study mappings defined by the stochastic integrals with respect to L\'evy processes
related to $K_{\al}(\rd)$.

Let $\al\in\R$ and 
$$
\ep _{\al}(u) = 
\begin{cases}
\int_u^{1}x^{-\al-1}dx, & \,\, 0<u<1, \\
0, & \,\, u\ge 1.
\end{cases}
$$
Then, when $\al \neq 0$,
$$
\ep _{\al }(u) = 
\begin{cases}
\al^{-1}(u^{-\al}-1), & \,\, 0<u<1, \\
0, & \,\, u\ge 1,
\end{cases}
$$
and when $\al =0$,
$$
\ep _{0 }(u) = 
\begin{cases}
\log u^{-1}, & \,\, 0<u<1, \\
0, & \,\, u\ge 1.
\end{cases}
$$
Let $\ep ^*_{\al}(t)$ be the inverse function of $\ep_{\al}(u)$, that is, $t=\ep_{\al}(u)$
if and only if $u=\ep_{\al}^*(t)$.
Note that 
$$
\ep_{\al}(0) =
\begin{cases}
(-\al)^{-1}, &\q  \al < 0,\\
\infty , & \q \al \ge 0.
\end{cases}
$$
Then, when $\al \neq 0$,
$$
\ep ^*_{\al }(t) = 
\begin{cases}
(1+\al t)^{-1/\al}, & \,\, 0<t < \ep_{\al}(0), \\
0, & \,\, t\ge \ep_{\al}(0),
\end{cases}
$$
and when $\al =0$,
$$
\ep ^*_{0 }(t) = e^{-t},\,\, t>0.
$$
Let $\{ X_t^{(\mu)}\}$ be the L\'evy process on $\rd$ with the distribuiton  $\mu\in I(\rd)$ at $t=1$.

\begin{defn}
Let $\al \in \R$.
We define mappings
$
\Phi _{\al} : 
{\frak D}(\Phi_{\al})  \,\,\to\,\, I(\rd)
$
by
$$
\Phi _{\al} (\mu) =\law\left (\int _0^{\ep_{\al}(0)} \ep_{\al}^*(t)dX_t^{(\mu)}\right ),
$$
where $\frak D (\Phi _{\al})$ is the domain of the mapping $\Phi _{\al}$.
\end{defn}

\begin{rem}
Let $-\infty<\beta<\alpha<\infty$. As in Sato (2006b) write the mapping
 as 
$$
\Phi_{\beta,\alpha}:= \law \left (\int_0^\infty f_{\beta,\alpha}(s)X_s^{(\mu)}\right ),
$$
where $f_{\beta,\alpha}(s)$ is the inverse function of 
$$
s=\left(
\Gamma(\alpha-\beta)
\right)^{-1}\int_t^1 (1-u)^{\alpha-\beta-1}u^{-\alpha-1}du.
$$
Our mappings in this paper $\Phi_\alpha$ are special cases of $\Phi_{\bt , \al}$ with $\bt =\al -1$.
Sato (2006b) discussed the domains of  $\Phi_{\bt , \al }$, but not the ranges of them,
and commented that description of the range of $\Phi_{\bt,\al}$ is to be made.
Our concern here is their ranges, although not for general $\bt < \al $,
because our motivation of this study started with the classes $K_{\al}(\rd)$.
\end{rem}
Regarding the domains of $\Phi _{\al}$,  we have the following result from Theorem 2.4 of Sato (2006b). 
\begin{prop}{\rm (Domains of $\Phi _{\al}$)}\\
$(i)$ When $\al <0$,  $\frak  D (\Phi_{\al})=I(\rd)$.\\
$(ii)$ When $\al =0$, $\frak D (\Phi_{\al})=I_{\log}(\rd)$.\\
$(iii)$ When $0<\al <1$, $\frak D (\Phi_{\al})=\{ \mu\in I(\rd) :
 \int_{\rd}|x|^\alpha \mu(dx)<\infty\} =: I_{\al}(\rd)$.\\
$(iv)$ When $\al =1$, $\frak D (\Phi_{1})= \{\mu\in I(\rd) : \int_{\rd}|x|\mu(dx)<\infty,$\\
\hskip 10mm $\lim_{T\to\infty}\int_{1}^Tt^{-1}dt
\int_{|x|>t}x\nu(dx)\,\,\text{exists in}\,\,\rd,$
 $\int_{\rd}x\mu (dx)=0 \}= :I_1^*(\rd)$.\\
$(v)$ When $1<\al <2$, $\frak D (\Phi_{\al})=\{ \mu\in I(\rd) :
 \int_{\rd}|x|^\alpha \mu(dx)<\infty, \int_{\rd}x\mu(dx)=0\} $\\
\hskip 10mm $=: I_\alpha^0(\rd)$.\\ 
$(vi)$ When $\al\ge 2$, $\frak D (\Phi_{\al})=\{\delta_0\}$, where $\delta _0$ 
is the distribution with the total mass at $0$.\\
\end{prop}

\n
Note that when $\al <0$, the interval of the integral is finite, so the stochastic integral exists
for any $\mu\in I(\rd)$ by a result in Sato (2006a).
Because of (vi) above, we are only interested in the case $\al<2$.
So, from now on, we assume that $\al<2$.

\begin{rem}
O'Connor (1979) mentioned the definition of $\Phi_{\al}, -1 < \al <2$, and stated without proofs that 
$\frak D(\Phi_{\al}) = I_{\al}(\R^1), 0<\al<1$, and $\frak D(\Phi_{\al}) = I_{\al}^0(\R^1), 1<\al<2$,
but he did not mention the case $\al=1$.
Actually, as we will see, the case $\al=1$ is the most difficult case to handle.
\end{rem}

\begin{rem} (Ranges)
We know
$$
\Phi_0(I_{\log}(\rd) ) = L(\rd)
\qquad \text{(Wolfe (1982) and others).}
$$
In Jurek (1985), it is shown that
$$
U(\rd) = \left\{ \law\left (\int_0^1 t dX_t^{(\mu)}\right ), \,\,\mu\in I(\rd)\right \}.
$$
But this is trivially the same as $\Phi_{-1}(I(\rd))$.
\end{rem}

In the following denote the mapped distribution by
$\widetilde\mu=\Phi_\alpha(\mu)=\widetilde \mu_{(\widetilde A,\widetilde\nu,\widetilde\gamma)}$
with polar decomposition $(\wt\ld, \wt\nu_{\xi})$.
We want to prove
\begin{thm}
\label{range-m0}
The ranges of the mapping $\Phi_{\al}$ are, \\
$(i)$ when $\al <0,\Phi_{\al}(I(\rd )) = K_{\al}(\rd)$,\\
$(ii)$ when $\al =0,\Phi_{0}(I_{\log}(\rd )) = K_{0}(\rd)$,\\
$(iii)$ when $ 0<\al<1, \Phi_{\al}(I_{\al}(\rd ))=K_{\al}(\rd)$,\\
$(iv)$ when $\al =1, \Phi_{\al}(I_{1}^*(\rd ))=$\\
\hskip 10mm $\{\widetilde \mu \in K_1(\rd) \,\,\text{with}\,\,\wt\nu_{\xi}(dr) = r^{-\al-1}\wt\ell_{\xi}(r)dr
\,\,\text{such that}$\\
\hskip 15mm $\lim_{\varepsilon \downarrow 0}
\int_\varepsilon^1t dt \int_S \xi\wt\ld(d\xi) \sek
 \frac{r^2}{1+t^2 r^2}d\wt\ell_{\xi}(r) \,\,\text{exists in $\rd$ and
 equals $\widetilde\gamma$}\}$,\\
$(v)$ when $1<\al<2, \Phi_{\al}(I_{\al}^0(\rd ))=K_{\al}(\rd)\cap \{\widetilde\mu\in 
I(\rd)\,\,\text{such that}\,\, \int_{\rd}x\widetilde\mu(dx)=0\}$.
\end{thm}

Although (ii) is known, we have written it just for the completeness of the theorem.
We give the proof of Theorem 4.7 in the next section.

We end this section with mentioning the continuity of $\Phi _{\al}(\mu)$ in $\al$ near $0$ 
from below for each fixed 
$\mu\in I_{\log}(\rd)$.
(The continuity in $\al \in [-1,0)$ for fixed $\mu\in I(\rd)$ is trivial.)

\begin{rem}
Now, let $\al$ tend to $0$ from below.
As to the interval of the integral, we have
$$
\int_0^{-1/\al} \to \sek\q\text{as}\,\, \al \uparrow 0
$$
and as to the integrand, we have
$$
(1+\al t)^{-1/\al} \to e^{-t}\q\text{as}\,\, \al \uparrow 0.
$$
So, the question is whether $\lim_{\al \uparrow 0} \Phi_{\al} (\mu) 
= \Phi_0(\mu), \mu\in I_{\log}(\rd)$, holds or not.
But, this is true, if we apply the dominated convergence theorem to the cumulants of $\Phi_{\al}(\mu)$.
\end{rem}

This remark explains why our expression (\ref{eq:law-alpha}) is more natural, when we consider 
the case $\al =0$ as mentioned in Remark 2.1.

\vskip 10mm

\section{Proof of Theorem 4.6}

Suppose $\mu=\mu_{(A,\nu,\gamma)}\in \frak D(\Phi_\alpha), -\infty<\alpha<2$. 
Then the mapped distribution
$\widetilde\mu=\Phi_\alpha(\mu)=\wt\mu_{(\widetilde A,\widetilde\nu,\widetilde \gamma)}$ satisfies 
\begin{equation}
\widetilde A=(2-\alpha)^{-1}A,
\label{eq:A-mapped}
\end{equation}
\begin{equation}
\widetilde \nu(B) = \int_0^1 \nu(s^{-1}B)s^{-\alpha-1}ds , \label{eq:nu-mapped}
\end{equation}
\begin{align}
\widetilde\gamma 
&= \lim_{\epsilon \downarrow 0}\int^1_\epsilon t^{-\alpha}dt \left(\gamma
-\int_{\mathbb R^d}x\left(\frac{1}{1+|x|^2}-\frac{1}{1+t^2|x|^2}\right)\nu(dx)\right)\label{eq:gamma-mapped} \\
&= \lim_{T\to\varepsilon_\alpha(0)} \int_0^T \varepsilon_\alpha^\ast (s)ds \left(\gamma-\int_{\mathbb R^d}x
\left(\frac{1}{1+|x|^2}-\frac{1}{1+|\varepsilon_\alpha^\ast(s)x|^2}\right)\nu(dx) \right).\nonumber
\end{align}
The derivation of $\widetilde A$ is that
\begin{align*}
\widetilde A &= \left (\int_0^{\varepsilon_\alpha(0)}\varepsilon_{\alpha}^\ast(t)^2dt \right ) A
= \left (\int_1^0 s^2 d\varepsilon_\alpha(s)\right )A
= (2-\alpha)^{-1}A.
\end{align*}
(5.2) is shown as follows.
By using 
Proposition 2.6 of Sato (2006b), we have
\begin{align*}
\wt\nu(B) 
& = \int _0^{\varepsilon_\alpha(0)} dt \int_{\rd} 1_B(x\ep_{\al}^*(t))\nu(dx)\\
& = \int_0^1 (-d\ep _{\al}(s))\int_{\rd}1_B(xs)\nu(dx)\\
& = \int_0^1 s^{-\al-1}ds \int_{\rd} 1_{s^{-1}B}(x)\nu(dx)\\
& = \int_0^1 \nu (s^{-1}B)s^{-\al-1}ds.
\end{align*}
{Similarly, by the change of variables $t\to\varepsilon_\alpha^\ast(s)$
we obtain two representations for $\tilde \gamma$.}
We sometimes use the zero mean condition, 
\begin{equation}
\gamma=-\int_{\mathbb{R}^d}\frac{x|x|^2}{1+|x|^2}\nu(dx).
\label{eq:cond-mean-zero}
\end{equation}
We need the following lemma. Denote 
\begin{align*}
\log^\ast x:= 
\left\{
\begin{array}{ll}
1 & \mbox{if}\quad 0<x\le 1,  \\
\log x &  \mbox{if}\quad x>1.
\end{array}
\right.
\end{align*}

\begin{lem}
\label{lem:eq:relation-widetilde-nu}
Let $-\infty<\alpha<2$ and let $\wt \nu$ be a L\'evy measure.
Then there exists a L\'evy measure $\nu$ satisfying $(5.2)$ such that
\begin{eqnarray}
\label{eq:relation-widetilde-nu}
\begin{cases}
\int_{\mathbb R^d}(|x|^2 \land 1)\nu(dx)<\infty, & \text{when}\ \alpha<0,\\
\int_{\mathbb R^d}(|x|^2 \land 1)\log^\ast|x| \nu(dx)<\infty, & \text{when}\ \alpha=0,\\
\int_{\mathbb R^d}(|x|^2 \land |x|^\alpha)\nu(dx)<\infty, & \text{when}\ 0<\alpha<2\\
\end{cases}
\end{eqnarray}
if and only if 
$\widetilde \nu$ is represented as 
\begin{equation}
\label{eq:range-infty-1}
\widetilde\nu (B) = \int_{S}\widetilde\lambda(d\xi)\int_0^\infty
 1_{B}(u\xi)u^{-\alpha-1}
\widetilde \ell_{\xi}(u)du, \quad B\in\mathcal{B}(\mathbb R^d),
\end{equation}
where $\widetilde \lambda$ is a measure on $S$ and $\widetilde
 \ell_{\xi}(u)$ is a function measurable in $\xi$ and for $\widetilde
 \lambda\text{-{\rm a.e.}}\,\xi.$ nonincreasing in $u\in(0,\infty)$, not identically zero
 and $\lim_{u\to\infty} \widetilde \ell_\xi(u)=0$.
\end{lem}

This lemma follows from similar arguments as those used in Lemma 4.4 in Sato (2006b).

\n
{\it Proof of Lemma \ref{lem:eq:relation-widetilde-nu}. }
We first show the \lq\lq only if''part.
Assume that the L\'evy measure $\nu$ satisfies (5.2) and  (\ref{eq:relation-widetilde-nu}). 
The polar decomposition gives us
\begin{eqnarray}
\int_{S}\lambda(d\xi)\int_0^\infty(r^2 \land 1)\nu_{\xi}(dr)<\infty,&& \
 \text{when}\ \alpha\le 0  \label{cond1-polar-decomp} \\
\int_{S}\lambda(d\xi)\int_0^\infty(r^2 \land
 r^\alpha)\nu_{\xi}(dr)<\infty, && \
 \text{when}\ \alpha> 0. \label{cond2-polar-decomp}
\end{eqnarray}
Then we have for $B\in\mathcal{B}(\mathbb R^d)$
\begin{align*}
\widetilde \nu(B) &= \int_0^1 \nu(s^{-1}B)s^{-\alpha-1}ds \\
&= \int_0^1 \int_{S} \lambda(d\xi) \int_0^\infty\nu_\xi(dr) 1_{s^{-1}B}(r\xi)s^{-\alpha-1}ds \\
&= \int_{S}\lambda(d\xi) \int_0^\infty r^\alpha \nu_\xi(dr)\int_0^r1_B(u\xi) u^{-\alpha-1}du\\ 
&=: \int_{S}\lambda(d\xi)\int_0^\infty 1_B(u\xi)u^{-\alpha-1}\widetilde \ell_{\xi}(u)du,
\end{align*}
where 
\begin{equation}
\widetilde \ell_{\xi}(u)=\int_u^\infty r^\alpha \nu_{\xi}(dr).
\end{equation}
Therefore $\widetilde \ell_{\xi}(u)$ is measurable in $\xi$, and for
$\lambda$-a.e.$\,\xi.$ nonincreasing in $u$, and $\lim_{u\to \infty}
\widetilde \ell_\xi(u)=0$ from (\ref{cond1-polar-decomp}) and
(\ref{cond2-polar-decomp}).

Conversely, suppose that $\widetilde \nu$ satisfies
(\ref{eq:range-infty-1}). Let $\widetilde\ell_\xi(u+)$ be the right-continuous function defined by 
$\lim_{t\uparrow u}\widetilde \ell_\xi(t)=\widetilde
\ell_\xi (u+)$. Then since $-\widetilde \ell_\xi(u+)$ is a right-continuous
increasing function, there exists a Lebesgue-Stieltjes measure $\widetilde Q_\xi$ on $(0,\infty)$ satisfying 
$$\widetilde Q_\xi((r,s])=-\widetilde \ell_\xi(s+) + \widetilde \ell_\xi(r+)$$ 
and put
$$
\nu_\xi(dr)=r^{-\alpha}\widetilde Q_\xi(dr).
$$
Furthermore, define
$$
\nu(B)=\int_{S}\widetilde \lambda(d\xi)\int_0^\infty 1_B(r\xi)\nu_\xi(dr).
$$
Let $\lambda=\widetilde \lambda$. 
Then for the case $\alpha<0$ we have
\begin{align*}
\int_0^\infty (|x|^2 & \land 1)\nu(dx) = \int_{S}\lambda(d\xi)\int_0^\infty(r^2\land 1) \nu_\xi(dr)\\
&= \int_{S}\lambda(d\xi) \left(\int_0^1r^{2-\alpha}\widetilde Q_\xi(dr)+\int_1^\infty
r^{-\alpha}\widetilde Q_\xi(dr)\right). 
\end{align*}
Since $\widetilde \nu$ is a L\'evy measure, we have $\int_{S}\lambda(d\xi)\int_0^\infty
(r^2\land 1)r^{-\alpha-1}\widetilde \ell_\xi(r+)dr<\infty$.
Note that
\begin{align*}
0< \int_{S}\lambda(d\xi)&\int_0^1r^{1-\alpha}\widetilde \ell_\xi(r+)dr = \int_{S}
\lambda(d\xi)\int_0^1r^{1-\alpha}\int_r^\infty \widetilde Q_\xi(dx)dr \\
&= \int_{S}\lambda(d\xi)\int_0^\infty \widetilde Q_\xi(dx)\int_0^{x\land 1} r^{1-\alpha}dr\\
&= ({2-\alpha})^{-1}\int_{S}\lambda(d\xi) \int_0^1x^{2-\alpha}\widetilde
 Q_\xi(dx)+({2-\alpha})^{-1}\int_{S}\lambda(d\xi){\widetilde \ell_\xi(1+)}<\infty
\end{align*}
and 
\begin{align*}
0<\int_{S}\lambda(d\xi) \int_1^\infty r^{-\alpha-1}&\widetilde \ell_\xi(r+)dr = \int_{S}\lambda(d\xi)\int_1^\infty
r^{-\alpha-1}\int_r^\infty \widetilde Q_\xi(dx)dr \\
&= \int_{S}\lambda(d\xi) \int_1^\infty \widetilde Q_\xi(dx)\int_1^x r^{-\alpha-1}dr \\
&= {\alpha}^{-1}\int_{S}\lambda(d\xi)\int_1^\infty ({1-x^{-\alpha}})\widetilde Q_\xi(dx) \\
&= {\alpha}^{-1}\int_{S}\lambda(d\xi) {\widetilde \ell_\xi(1+)}- {\alpha}^{-1}\int_{S}
\lambda(d\xi)\int_1^\infty x^{-\alpha} \widetilde Q_\xi(dx)<\infty.
\end{align*}
From the first inequality, $\int_{S}
\lambda(d\xi)\widetilde \ell_\xi(1+)>0 $ is finite and we see that 
$$
0< \int_{S}\lambda(d\xi) \int_0^1x^{2-\alpha}\widetilde
 Q_\xi(dx)<\infty\quad \text{and} \quad 0<\int_{S}\lambda(d\xi)\int_1^\infty x^{-\alpha} \widetilde 
Q_\xi(dx)<\infty,
$$
which imply {(\ref{eq:relation-widetilde-nu}).} 
{For the remaining cases $\alpha=0$ and $0<\alpha<2$, similar logic as in
the case $\alpha<0$ works and we concludes
(\ref{eq:relation-widetilde-nu}).} 
\hfill $\Box$

\vskip 3mm
\noindent
{\it Proof of Theorem 4.6.}
As in Sato (2006a), 
we use the notation $C^+_{\#}$ for the class of nonnegative bounded continuous
functions on $\mathbb R^d$ vanishing on a neighborhood of the origin.\\
(i), (ii) and (iii) $(-\infty<\alpha<1)$ \\
Suppose that $\widetilde \mu\in \Phi_{\al}(\mathfrak D(\Phi_\alpha))$ and $\wt\mu = \Phi_\al(\mu), 
\mu=\mu_{(A,\nu,\gm)}$.
When $\nu\ne 0$, since  $\mu\in\mathfrak D (\Phi_{\al})$, (\ref{eq:range-infty-1}) holds by
Lemma \ref{lem:eq:relation-widetilde-nu} so that $\wt\mu\in K_{\al}(\rd)$.

Conversely, suppose $\wt\mu = \wt\mu_{(\wt A, \wt\nu, \wt\gm)}\in K_{\al}(\rd)$.
If $\widetilde\mu$ is Gaussian then putting $A=(2-\alpha) \widetilde A, \nu=0$, and
$\gamma=(1-\alpha)\widetilde \gamma$, we have
$\mu=\mu_{(A,\nu,\gamma)}\in \mathfrak D(\Phi_\alpha)$ and $\widetilde\mu=\Phi_\alpha(\mu)$. 
If $\widetilde\mu$ is non-Gaussian, then we have (\ref{eq:relation-widetilde-nu}) 
by Lemma \ref{lem:eq:relation-widetilde-nu}. 
We put $A=(2-\alpha) \widetilde A$ and 
$$
\gamma=(1-\alpha)\left(
\widetilde\gamma+\int_0^{\varepsilon_\alpha(0)}
\varepsilon_\alpha^\ast(t)dt\int_{\mathbb{R}^d}x\left(
\frac{1}{1+|x|^2}-\frac{1}{1+|\varepsilon_\alpha^\ast(t)x|^2}\right)\nu(dx)\right).
$$
Although the parametrization of $\alpha$ is different, the argument
similar to the proof of (2.35) in 
Sato (2006b) works and
it follows from (\ref{eq:relation-widetilde-nu}) that
$$
\int_0^{\varepsilon_\alpha(0)} \varepsilon_\alpha^\ast(t) dt \int_{\mathbb{R}^d}|x|\left|
\frac{1}{1+|x|^2}-\frac{1}{1+|\varepsilon_\alpha^\ast(t) x|^2}\right|\nu(dx)<\infty.
$$
Thus $\mu=\mu_{(A,\mu,\gamma)}\in
\mathfrak D(\Phi_\alpha)$ and $\Phi_\alpha(\mu)=\widetilde\mu$. \\
(iv) $(\alpha=1)$ \\
Suppose that $\widetilde \mu =\wt\mu _{(\wt A,\wt \nu,\wt\gm)} = \Phi (\mu) \in\Phi_1(\mathfrak D(\Phi_1))$ and
$\mu=\mu_{(A,\nu,\gamma)}\in \mathfrak D(\Phi_1)$. First assume that
$\widetilde \mu$ is Gaussian. 
Then for given $\varphi \in C_{\#}^+$,  
$
0=\int_0^1 \int_{\mathbb R^d}\varphi(sx)s^{-2}\nu(dx)ds,
$
which implies $0=s^{-2}\int_{\mathbb{R}^d}\varphi(sx)\nu(dx)$ a.e. 
Since by the dominated convergence theorem $s^{-2}\int_{\mathbb
R^d}\varphi(sx)\nu(dx)$ is continuous in $s$, letting $s=1$, we have $\nu=0$. 
Furthermore, from Proposition 4.3 (iv) with (\ref{eq:cond-mean-zero}) 
$\gamma=0$ and hence $\widetilde \gamma=0$. 
When $\widetilde \mu$ is non-Gaussian, $\nu$  satisfies (\ref{eq:relation-widetilde-nu}) with
$\alpha=1$, and (\ref{eq:gamma-mapped}) and (\ref{eq:cond-mean-zero}) imply that
\begin{equation}
\label{eq:cond-conve-gamma-mapped}
-\lim_{\varepsilon \downarrow 0}
\int_\varepsilon^1 t dt \int_{\mathbb R^d}\frac{x|x|^2}{1+t^2|x|^2}\nu(dx)
\end{equation}
exists in $\mathbb R^d$ and equals $\widetilde \gamma$.
Thus, $\wt\mu \in K_1(\rd)\cap \{\mu \in I(\rd) \,\,\text{such that}$
$\linebreak$ $-\lim_{\varepsilon \downarrow 0}
\int_\varepsilon^1t dt \int_{\mathbb R^d}\frac{x|x|^2}{1+t^2|x|^2}\nu(dx)\,\,\text{exists in $\rd$}\}$.

We show the converse. 
Suppose $\wt\mu =\wt \mu_{(\wt A, \wt \nu, \wt\gm)} \in K_1(\rd)\cap \{\mu \in I(\rd) \,\,\text{such that}$
$\linebreak$ $-\lim_{\varepsilon \downarrow 0}
\int_\varepsilon^1t dt \int_{\mathbb R^d}\frac{x|x|^2}{1+t^2|x|^2}\nu(dx)\,\,\text{exists in $\rd$}\}$.
If $\widetilde \mu$ is centered Gaussian, then
$\widetilde \mu\in\Phi_1(\mathfrak D(\Phi_1))$ from Proposition 4.3. 
If $\widetilde \mu$ is non-Gaussian and satisfies (\ref{eq:range-infty-1})
and (\ref{eq:cond-conve-gamma-mapped}), then by Lemma
\ref{lem:eq:relation-widetilde-nu} a measure $\nu$ exists and satisfies
(\ref{eq:nu-mapped}) and 
(\ref{eq:relation-widetilde-nu}) with $\alpha=1$. 
Let $\gamma=-\int_{\mathbb R^d}\frac{x|x|^2}{1+|x|^2}\nu(dx)$ and $A=\widetilde A$.
It follows from the existence of (\ref{eq:cond-conve-gamma-mapped}) and $\int_{|x|>1}|x|\nu(dx)<\infty$ that
$$
\lim_{T\to\infty}\int_{t_0}^T t^{-1} dt\int_{|x|>t}x\nu(dx)<\infty
$$
as in the proof of Theorem 2.8 of Sato (2006b).
Thus $\mu\in\mathfrak D(\Phi_1)$. Furthermore
(\ref{eq:cond-conve-gamma-mapped}) implies 
\begin{align*}
\widetilde \gamma &= \lim_{\varepsilon \downarrow 0}\int_\varepsilon^1 t^{-1}dt \left(
-\int_{\mathbb R^d}\frac{x|x|^2}{1+|x|^2}\nu(dx) +\int_{\mathbb R^d}x\left(\frac{1}{1+|tx|^2}
-\frac{1}{1+|x|^2}\right)\nu(dx)\right),
\end{align*}
which equals the right-hand side of (\ref{eq:gamma-mapped}). 
Therefore $\Phi_1(\mu)=\widetilde\mu$ and $\widetilde \mu\in \Phi_1(\mathfrak D(\Phi_1))$. \\
(v) $(1<\alpha<2)$ \\
Assume that $\widetilde \mu=\Phi_\alpha (\mu)$ with some
$\mu=\mu_{(A,\nu,\gamma)}\in\mathfrak D(\Phi_\alpha)$. The Gaussian case
is the same as that of the proof for (ii). 
If $\widetilde \mu$ is non-Gaussian, then it follows from Lemma
\ref{lem:eq:relation-widetilde-nu} that there exists $\widetilde \nu$
satisfying (\ref{eq:range-infty-1}). Since $\mu\in \mathfrak D(\Phi_\alpha),\nu$ and $\gamma$
satisfy $\int_{|x|>1}|x|^\alpha\nu(dx)<\infty$ and
(\ref{eq:cond-mean-zero}), respectively. 
Then as in the proof of Theorem 2.4 (iii) of Sato (2006b), $\widetilde\gamma$ exists and equals to 
\begin{align}
\widetilde \gamma = -\int_0^\infty \varepsilon_\alpha^\ast(t)dt\int_{\mathbb
R^d}\frac{x|\varepsilon_\alpha^\ast(t)x|^2}{1+|\varepsilon_\alpha^\ast(t)x|^2}\nu(dx)\noindent  
=-\int_{\mathbb R^d}\frac{x|x|^2}{1+|x|^2}\widetilde \nu(dx),
\label{eq:mean-zero-mapped}
\end{align}
which is 
\begin{equation}
\label{eq:mean-zero}
\int_{\mathbb R^d} x\widetilde \mu(dx)=0. 
\end{equation}
Hence $\wt\mu \in K_{\al}(\rd)\cap \{\mu\in I(\rd) : \int_{\rd}x\mu(dx)=0\}$.

We show the converse. 
Suppose $\wt\mu = \wt\mu_{(\wt A, \wt\nu, \wt\gm)} \in K_{\al}(\rd)\cap \{\mu\in I(\rd) : \int_{\rd}x\mu(dx)=0\}$.
The Gaussian case is obvious. Suppose $\widetilde
\mu$ be non-Gaussian. 
Due to Lemma
\ref{lem:eq:relation-widetilde-nu} a measure $\nu$
with $\nu(\{0\})=0$ exists and satisfies (\ref{eq:nu-mapped}) and
(\ref{eq:relation-widetilde-nu}). 
It follows from (\ref{eq:nu-mapped}) that 
\begin{align*}
\int_{\mathbb R^d}\frac{|x|^3}{1+|x|^2}&\widetilde \nu(dx) = 
\int_0^1 t^{2-\alpha}dt \int_{\mathbb R^d}
\frac{|x|^3}{1+t^2|x|^2}\nu(dx) \\
&\le \int_{|x|\le 1}|x|^3\nu(dx)\int_0^1t^{2-\alpha}dt +\int_{|x|>1}|x|^3 \nu(dx) \int_0^{1/|x|}t^{2-\alpha}dt \\
&\hskip 20mm +\int_{|x|>1}|x|\nu(dx)\int_{1/|x|}^1t^{-\alpha} dt \\
&=(3-\alpha)^{-1}\int_{|x|\le1}|x|^3\nu(dx)+({3-\alpha})^{-1}\int_{|x|>1}|x|^\alpha \nu(dx)\\
&\hskip 20mm +({1-\alpha})^{-1}\int_{|x|>1}\left(|x|-|x|^\alpha \right)\nu(dx)<\infty.
\end{align*}
Hence we have $\int_{|x|>1}|x|\widetilde \nu(dx)<\infty$ 
which is equivalent to $\int_{\mathbb R^d}|x|\widetilde\mu(dx)<\infty$ and
(\ref{eq:mean-zero-mapped}) holds. Let $\gamma=-\int_{\mathbb
R^d}\frac{x|x|^2}{1+|x|^2}\nu(dx),\
A=(2-\alpha) \widetilde A$ and
$\mu=\mu_{(A,\nu,\gamma)}$. Then $\mu\in \mathfrak D(\Phi_\alpha)$ by
Proposition 4.4 (v). Further 
\begin{eqnarray*}
\int_0^\infty \varepsilon_\alpha^\ast(t)dt\left(
\gamma+\int_{\mathbb R^d}x \left(
\frac{1}{1+|\varepsilon_\alpha^\ast(t)x|^2}-\frac{1}{1+|x|^2}
\right)\nu(dx)
\right)
=-\int_{\mathbb R^d}\frac{x|x|^2}{1+|x|^2}\widetilde \nu(dx),
\end{eqnarray*}
which equals $\widetilde \gamma$. Hence (\ref{eq:gamma-mapped}) is true
and $\Phi_\alpha(\mu)=\widetilde \mu$, namely $\widetilde \mu \in\Phi_{\al}(\mathfrak D(\Phi_\alpha))$. 
\qed

\vskip 10mm

\section{Nested subclasses of $\Phi_{\al}(\frak D(\Phi_{\al}))$}

$\Phi_{\al}$-mapping allows us to construct 
nested subclasses of $\Phi_{\al}(\frak D(\Phi_{\al}))$ denoted by $\Phi
_{\al}^{m+1},\ m=1,2,\ldots.$ This is the topic in this section. We will
see the domains $\frak D (\Phi _{\al}^{m+1})$ in subsection \ref{sebsec:domain:nested:subclass}
and characterize the ranges $\Phi_{\al}^{m+1}(\frak
D(\Phi_{\al}^{m+1}))$ by both stochastic integral representations and
the L\'evy-Khintchine triplet, which are respectively given in subsections
\ref{subsec:range1:nested:subclass} and
\ref{subsec:range2:nested:subclass}. 

\vskip 3mm
\subsection{ Domains of $\Phi _{\al}^{m+1}$}
\label{sebsec:domain:nested:subclass}

\begin{thm} Let $m=1,2,...$
\label{thm:integral-representation-domain-K}\\
$(i)$  When $\al <0$,
$
\frak D (\Phi _{\al}^{m+1}) = I(\rd).
$\\
$(ii)$ When $\al =0$,
$$
\frak D (\Phi _{0}^{m+1}) = \left\{ \mu\in I(\rd)
 : \int_{\rd} (\log |x|)^{m+1}\mu(dx)< \infty\right \} =: I_{\log^{m+1}}(\rd).
$$
$(iii)$ When $0<\al < 1$,
$$
\frak D (\Phi _{\al}^{m+1}) = \left\{
\mu\in I(\rd) :\int_{\rd}|x|^\alpha \left(\log|x|\right)^{m}\mu(dx)<\infty
\right\} =:I_{\al, \log^m}(\rd).
$$
$(iv)$ When $\alpha=1$,
\begin{align*}
\frak D (\Phi _{1}^{m+1}) &= \left\{
\mu\in I(\mathbb R^d) :\int_{\mathbb R^d}|x|
 \left(\log|x|\right)^{m}\mu(dx)<\infty,\ \int_{\mathbb R^d} x\mu(dx)=0,
				   \right. \\
& \left. \hskip -10mm \lim_{T\to\infty}\int_{t_0}^Tt^{-1}dt
\int_{|x|>t}x(\log(|x|/t))^m\nu(dx)\,\,\text{exists in}\,\,\mathbb R^d
\right\} =: I_{1,\log^m}^*(\rd).
\end{align*}
$(v)$ When $1<\al< 2$,
\begin{align*}
\frak D (\Phi _{\al}^{m+1}) &= \left\{
\mu\in I(\rd) :\int_{\rd}|x|^\alpha
 \left(\log|x|\right)^{m}\mu(dx)<\infty,\ \int_{\rd} x\mu(dx)=0
\right\} \\
&=:I^0_{\al, \log^m}(\rd).
\end{align*}
\end{thm}

\n
{\it Proof of Theorem 6.1.}
Since when $ \al <0$, the integral for $\Phi_{\al}^{m+1}(\mu)$ is not improper integral, it is easy 
to see that $\frak D(\Phi_{\al}^{m+1})=I(\rd)$, (see Sato (2006a)).
 When $\al=0$, Jurek (1985) determined $\frak D(\Phi_0^{m+1})$ 
as above.  

We are now going to prove (iii), (iv) and (v).
First, note that
\begin{equation}
\label{equation:integral-rep-domain-K-rem}
\int^1_{1/x} u^{-1}(\log ux)^m du = {(m+1)^{-1}}(\log x)^{m+1} \,\,
\text{for}\,\,m\ge0.
\end{equation}
Now, Theorem \ref{thm:integral-representation-domain-K} (iii), (iv) and
(v) are true for $m=0$ as seen in Proposition 4.3 (iii), (iv) and (v).  
Suppose that it is true for some $m>0$, as the induction hypothesis.  
Suppose $0<\al<1$.  Then
\begin{align*}
\frak D(\Phi^{m+2}_{\al})=\{ \mu \in \frak D(\Phi_{\al}) :& \int_{|x|>1}|x|^{\al}(\log |x|)^m \wt \nu(dx)<\infty, \\
&\text{where}\,\, \wt\nu \,\,\text{is the L\'evy measure of} \,\, \wt\mu=\Phi_{\al}(\mu) \}.
\end{align*}
Recall from (5.2) that
$$
\wt \nu(B)=\int_0^1 \nu(s^{-1}B)s^{-\al-1}ds.
$$
Thus, 
\begin{align*}
\int_{|x|>1} |x|^{\al} &(\log |x|)^m \wt \nu(dx) \\
=&\int_{|x|>1} |x|^{\al}(\log |x|)^m \int^1_0 \nu(s^{-1}dx)s^{-\al-1}ds \\
=&\int^1_0s^{-\al-1}ds \int_{|x|>1} |x|^{\al}(\log |x|)^m \nu(s^{-1}dx) \\
=&\int_{|y|>1}|y|^{\al}\nu(dy) \int^1_{1/|y|}s^{-1}(\log|sy|)^m ds.
\end{align*}
Then by (\ref{equation:integral-rep-domain-K-rem}),
$$
\int_{|x|>1}|x|^{\al}(\log |x|)^m \wt \nu(dx) <\infty
$$
if and only if
$$
\int_{|x|>1}|x|^{\al}(\log |x|)^{m+1} \nu(dx)<\infty,
$$
and we conclude that
$
\frak D(\Phi^{m+2}_\al)=I_{\al,\log^{m+1}}(\R^d).
$

When $1<\al<2$, there is no problem for the moment condition, and the condition, 
$\int_{\rd}x\mu(dx)=0,$ always holds.
Thus we get $\frak D(\Phi_{\al}^{m+2})=I_{\al, \log^{m+1}}^0(\rd)$. 

Finally we prove (iv). So, suppose $\al =1$.
Also suppose it is true for some
$m>0$. We have
\begin{align}
\frak D (\Phi _{1}^{m+2}) &= \{
\mu\in \frak D (\Phi _{1}) :\int_{|x|>1}|x|
 \left(\log|x|\right)^{m}\wt \nu(dx)<\infty,\ \int_{\mathbb R^d} x\mu(dx)=0,
				    \nonumber \\
& \lim_{T\to\infty}\int_{t_0}^Tt^{-1}dt
\int_{|x|>t}x(\log(|x|t^{-1}))^m\wt \nu(dx)\,\,\text{exists in}\,\,\mathbb R^d, 
\label{ineq:condition-domain-K-(iv)}\\
& \text{where}\,\, \wt\nu\,\,\text{is the L\'evy measure of} \,\,
 \Phi_{1}(\mu)\}. \nonumber
\end{align}
Since the moment condition can be given by the same way as for the case $1<\al<2$, in order to
reach the conclusion, it remains to show that
\begin{equation}
\label{equation:condition-m+1}
\lim_{T\to\infty}\int_{t_0}^Tt^{-1}dt
\int_{|x|>t}x(\log(|x|t^{-1}))^{m+1} \nu(dx)\,\,\text{exists in}\,\,\mathbb R^d.
\end{equation}
We have
\begin{align*}
\int_{|y|>t}y&(\log(|y|t^{-1}))^m\wt \nu(dy) \\
&= \int_{|y|>t}y(\log(|y|t^{-1}))^m \int_0^1 \nu(s^{-1}dy)s^{-2}ds \\
&= \int_0^1 s^{-1}ds \int_{|sx|>t} x(\log(|sx|t^{-1}))^m \nu(dx) \\
&= \int_{|x|>t}x \nu(dx) \int_{t/|x|}^1s^{-1} (\log(s|x|t^{-1}))^m ds\\
&= ({m+1})^{-1} \int_{|x|>t} x (\log(|x|t^{-1}))^{m+1} \nu(dx).
\end{align*}
Hence (\ref{ineq:condition-domain-K-(iv)}) is equivalent to (\ref{equation:condition-m+1}).
This completes the proof.
\qed

\vskip 3mm

\subsection{ Characterizations of  $\Phi_{\al}^{m+1}(\frak D(\Phi_{\al}^{m+1}))$
by stochastic integral representations}
\label{subsec:range1:nested:subclass}

\begin{thm}
\label{thm:integral-representation-K}
Let $\al <2$ and
$$
g_{\al,m}(s) = (m!)^{-1}s^{-\al-1}(\log s^{-1})^{m}1_{(0,1]}(s),
$$
$$
\ep_{\al,m}(u) = \int_u^{\infty}g_{\al,m}(s)ds = 
\begin{cases}
\int_u^{1}g_{\al, m}(s)ds,& \,\, 0<u<1, \\
0 ,& \,\, u\ge 1,
\end{cases}
$$
and let 
$\ep_{\al,m}^*(t)$ be the inverse function of $\ep_{\al,m}(x)$ such that  $ t=\ep_{\al,m}(u)$ if and only if
$u=\ep_{\al,m}^*(t)$.
Then
$\wt\mu\in \Phi_{\al}^{m+1}(\frak D(\Phi_{\al}^{m+1}))$ if and only if
\begin{equation}
\wt\mu = \law \left ( \int _0^{\ep_{\al,m}(0)} \ep^*_{\al,m}(t)dX_t^{(\mu)}\right ),
\,\,\text{for some}\,\,\mu\in \frak D (\Phi_{\al}^{m+1}),
\end{equation}
where when $ \al <0$, $\ep_{\al,m}(0)= (-\al)^{-(m+1)}$ and when $0\le \al <2$, $\ep_{\al,m}(0)=\infty$.
\end{thm}

\begin{rem}
When $ \al <0$,  $\ep_{\al,m}(0)= (-\al)^{-(m+1)}$ above is shown as follows.
\begin{align*}
\ep_{\al,m}(0) = \int_0^{1}g_{\al, m}(s)ds = (m!)^{-1}\sek e^{\al t}t^mdt = (-\al)^{-(m+1)}.
\end{align*}
\end{rem}
\begin{rem}
The following are known.

\n
(1) (Jurek (2004).)
When $\al =-1$, 
$\wt\mu\in \Phi_{-1}^{m+1}(\mathfrak D (\Phi_{-1}^{m+1}))$ if and only if
\begin{equation*}
\wt\mu = \law \left ( \int _0^{1} \tau^*_{m}(t)dX_t^{(\mu)}\right ) \,\,\text{for some}\,\,\mu\in I(\rd),
\end{equation*}
where $\tau_m(u) = \int_0^u g_{-1,m}(s)ds, 0<u\le 1$ and $\tau_m^*(t)$ is its  inverse.
However, by changing variable $t$ to $1-t$, we see that
$$
\law \left ( \int _0^{1} \tau^*_{m}(t)dX_t^{(\mu)}\right ) = 
\law \left ( \int _0^{1} \ep^*_{-1,m}(t)dX_t^{(\mu)}\right ).
$$

\n
(2) (Jurek (1983).)
When $\al =0$,
\begin{equation}
\ep_{0,m}^*(t) = e^{-((m+1)!t)^{(m+1)^{-1}}} .
\end{equation}
In our setting, we can get (6.4) as follows.
By a standard calculation, we see that
$$
\ep_{0,m}(u) = ((m+1)!)^{-1}(\log u^{-1})^{m+1}
$$
and thus (6.4) is given by taking the inverse function of $t= \ep_{0,m}(u)$.
\end{rem}

We now prove Theorem 6.2.

\n
{\it Proof of Theorem 6.2.}
(\lq\lq Only if" part.) Let $\wt\mu \in \Phi_{\alpha}^{m+1}(\mathfrak D (\Phi_{\alpha}^{m+1}))$.
Then $\wt\mu = \Phi_{\al}^{m+1}(\mu)$ for some $\mu \in \frak D(\Phi^{m+1}_{\al})$.
We regard
$\Phi_{\al}$ as a mapping from a L\'evy measure to a L\'evy measure.
Namely,
$$
\Phi_{\al}(\nu)(B) := \nu_{\Phi_{\al}(\mu)}(B) = \int_0^1\nu_{\mu}(s^{-1}B)s^{-\al-1}ds,
$$
where $\nu_{\mu}$ is the L\'evy measure of $\mu\in I(\rd)$.
We first show, for each L\'evy measure $\nu$,
\begin{align}
{\Phi_{\al}^{m+1}(\nu)}(B) & = (m!)^{-1}\int_0^1\nu(s^{-1}B)s^{-\al-1}(\log s^{-1})^{m} ds\\
& = \int_0^1\nu(s^{-1}B)g_{\al,m}(s)ds, \quad m\ge 1. \nonumber
\end{align}
We have
\begin{align*}
\Phi_{\al}^2(\nu)(B) 
& = \Phi_{\al}(\Phi_{\al}(\nu))(B)\\
& = \int_0^1 \Phi_{\al}(\nu)(s^{-\alpha-1}B)s^{-\al-1}ds\\
& = \int_0^1s^{-\alpha-1}ds\int_0^1\nu((ts)^{-1}B)t^{-\al-1}dt\\
& = \int_0^1\nu(u^{-1}B)u^{-\al-1}du \int_u^1s^{-1}ds\\
& = \int_0^1\nu(u^{-1}B)g_{\al,1}(u)du.
\end{align*}
Thus (6.3) is true for $m=1$.
Next suppose
$$
{\Phi_{\al}^m(\nu)}(B) = ((m-1)!)^{-1}\int_0^1\nu(s^{-1}B)s^{-\al-1}(\log s^{-1})^{m-1} ds.
$$
Then, 
\begin{align*}
\Phi_{\al}^{m+1}&(\nu)(B) = \int_0^1 \Phi_{\al}^m (\nu)(s^{-1}B)s^{-\al-1}ds\\
& = ((m-1)!)^{-1}\int_0^1s^{-\al-1}ds \int_0^1\nu((us)^{-1}B)u^{-\al-1}(\log u^{-1})^{m-1} du \\
& = ((m-1)!)^{-1}\int_0^1u^{-\al-1}(\log u^{-1})^{m-1} du \int_0^u\nu(t^{-1}B)(t/u)^{-\al-1}(u^{-1})dt\\
& = ((m-1)!)^{-1}\int_0^1\nu(t^{-1}B)t^{-\al-1}dt \int_t^1 u^{-1}(\log u^{-1})^{m-1}du\\
& = (m!)^{-1}\int_0^1\nu(t^{-1}B)t^{-\al-1}(\log t^{-1})^{m} dt\\
& = \int_0^1\nu(t^{-1}B)g_{\al,m}(t)dt.
\end{align*}
$A$ of $\mu_{(A,\nu,\gamma)}$ is determined by $\widetilde A$ defined in
(\ref{eq:m-A-mapped}). Regarding $\gamma$, when $\alpha<1$, it is given
by $\widetilde \gamma$ in (\ref{eq:m-gamma-mapped}) since
$\int^1_0 t^{-\alpha}(\log t^{-1})^m dt$ is integrable. 
(Here, (6.7) and (6.9) will be given in the next section.)
When $1\le
\alpha <2$, from zero mean condition of $\mathfrak D
(\Phi_{\alpha}^{m+1})$, we have 
$\gamma=-\int_{\mathbb{R}^d}\frac{x|x|^2}{1+|x|^2}\nu(dx).$
 Hence $\widetilde \mu$ has the representation (6.4).

In general, (see, e.g. the equation (1.8) of Barndorff-Nielsen and Maejima (2008)), if 
$$
\wt\nu(B) = \sek \nu(s^{-1}B)g(s)ds,
$$
for some nonnegative integrable function $g(u)$ assuring the convergence of the integral,
then $\wt\mu\in I(\rd)$ with the L\'evy measure $\wt\nu$ satisfies
$$
\wt\mu = \law\left ( \int_0^{\ep(0)}\ep^*(t)dX_t\right ),
$$
where $\ep^*(t)$ is the inverse function of $\ep(x) = \int_x^{\infty}g(s)ds$ and
$\{X_t\}$ is a L\'evy process with its L\'evy measure $\nu$.
This is (6.4).

(\lq\lq If" part.)
Let 
$$\wt\mu = \law \left(\int^{\ep_{\al, m}(0)}_{0} \ep^*_{\al,m}(t)dX_t^{(\mu)}\right)
\,\,\text{for some}\,\, \mu \in \frak D(\Phi_{\al}^{m+1}).$$
Then
\begin{align*}
\nu_{\wt\mu}(B) 
&= \int_0^{\ep_{\al, m}(0)}dt \int_0^\infty 1_B(x\ep^*_{\al,m}(t)) \nu_\xi(dx)\\
&= \int_0^1(-d\ep_{\al,m}(u)) \int_0^\infty 1_B(xu)\nu_\xi(dx)\\
&= \int_0^1 g_{\al,m}(u)du \int_0^\infty 1_{u^{-1}B}(x)\nu_\xi(dx)\\
&= \int_0^1\nu(u^{-1}B) g_{\al,m}(u)du, 
\end{align*}
where when $ \al <0$, $\ep_{\al,m}(0)= (-\al)^{-(m+1)}$ and when $0\le \al <2$, $\ep_{\al,m}(0)=\infty$. 
$\widetilde A$ and $\widetilde \gamma$ are given respectively by (\ref{eq:m-A-mapped}) and
(\ref{eq:m-gamma-mapped}). 
Thus, $\wt\mu \in \Phi_{\alpha}^{m+1}(\mathfrak D (\Phi_{\alpha}^{m+1})).$
The proof is completed.
\qed

\vskip 3mm

\subsection{{ Characterizations of $\Phi_{\al}^{m+1}(\frak D(\Phi_{\al}^{m+1}))$ 
by the L\'evy-Khintchine triplet}} 
\label{subsec:range2:nested:subclass}

We consider the range $\Phi_{\al}^{m+1}(\frak D(\Phi_{\al}^{m+1}))$ for $m=1,2,\ldots$
Let $-\infty <\al <2$ and suppose $\mu=\mu_{(A,\nu,\gamma)}\in \mathfrak D(\Phi_{\al}^{m+1})$.
Then the mapped distribution
$\widetilde\mu=\Phi_{\al}^{m+1}(\mu) = \mu_{(\widetilde A,\widetilde \nu, \widetilde \gamma)}$ satisfies 
\begin{equation}
\widetilde A=(2-\alpha)^{-(m+1)}A,
\label{eq:m-A-mapped}
\end{equation}
\begin{equation}
\widetilde \nu(B) = (m!)^{-1} \int_0^1 \nu(s^{-1}B)s^{-\alpha-1}(\log s^{-1})^m ds , \label{eq:m-nu-mapped}
\end{equation}
\begin{align}
\widetilde\gamma 
&= 
\lim_{\epsilon \downarrow 0}\ (m!)^{-1} \int^1_\epsilon t^{-\alpha}(\log t^{-1})^m dt\left(
\gamma-\int_{\mathbb R^d}x\left(\frac{1}{1+|x|^2}-\frac{1}{1+t^2|x|^2}\right)\nu(dx)
\right) \nonumber \\
& \label{eq:m-gamma-mapped} \\
&= 
\lim_{T\to\varepsilon_{\alpha,m}(0)} \int_0^T \varepsilon_{\alpha,m}^\ast
 (s)ds \left(\gamma-\int_{\mathbb R^d}x\left(\frac{1}{1+|x|^2}-
\frac{1}{1+|\varepsilon_{\alpha,m}^\ast(s)x|^2}
\right)\nu(dx) \right).\nonumber
\end{align}
For $\widetilde A$, we used
\begin{equation*}
 \int_0^{\varepsilon_{\alpha,m}(0)}  \varepsilon_{\alpha,m}^\ast(t)^2dt 
= \int_1^0 s^2 d\varepsilon_{\alpha,m}(s) 
= (m!)^{-1}\int_0^1 s^{1-\alpha}(\log s^{-1})^m ds \\
={(2-\alpha)^{-(m+1)}}.
\end{equation*}
For $\widetilde \gamma$, we use the same calculation as that for
$\Phi_{\al}(\frak D(\Phi_{\al}))$.

\begin{thm}
\label{thm:range-of-K-m-mapping}
Let $-\infty<\alpha<2$ and $m=1,2,\ldots$ Then $\widetilde\mu\in
\Phi_{\al}^{m+1}( \mathfrak D(\Phi_\alpha^{m+1}))$ if and only if one of the following
 conditions depending on $\alpha$ is satisfied. \\
$(i)$ $(-\infty<\alpha<1)$ \\
$\widetilde\mu$ is Gaussian, or $\widetilde\mu$ is non-Gaussian and 
\begin{equation}
\label{eq:m-range-infty-1}
\widetilde\nu (B) = \int_{S}\widetilde\lambda(d\xi)\int_0^\infty
 1_{B}(u\xi)u^{-\alpha-1}
\widetilde h^{(m)}_{\xi}(u)du, \quad B\in\mathcal{B}(\mathbb R^d).
\end{equation}
Here $\widetilde \lambda$ is a measure on $S$ and $\widetilde
 h^{(m)}_{\xi}(u)$ is a measurable function in $\xi$ such that satisfies 
\begin{equation}
\label{eq:m-h}
\widetilde h^{(m)}_\xi(u)=((m-1)!)^{-1} \int_u^\infty {x}^{-1}\left(\log (x/u)\right)^{m-1}\widetilde \ell_\xi(x)dx
\end{equation}
where $\widetilde \ell_\xi(u)$ is a function measurable in $\xi$ and for 
$\lambda-a.e.\xi.$ nonincreasing in $u\in(0,\infty)$, not identically zero
and $\lim_{u\to\infty} \widetilde \ell_\xi(u)=0$. \\
$(ii)$ $(\alpha=1)$ \\
$\widetilde\mu$ is centered Gaussian, or $\widetilde \mu$ is
 non-Gaussian and $\widetilde \mu$ satisfies $(\ref{eq:m-range-infty-1})$, $(\ref{eq:m-h})$ and 
\begin{equation}
\label{eq:cond-conve-gamma-m-mapped}
-\lim_{\varepsilon \downarrow 0}
\int_\varepsilon^1t \left(\log t^{-1}\right)^m dt \int_{\mathbb R^d}\frac{x|x|^2}{1+t^2|x|^2}\nu(dx)
\end{equation}
exists in $\mathbb R^d$ and equals $\widetilde \gamma$. Here
the measure $\nu$ satisfying $(\ref{eq:m-nu-mapped})$.\\
$(iii)$ $(1<\alpha<2)$ \\
$\widetilde\mu$ is centered Gaussian, or $\widetilde\mu$ is non-Gaussian
 and $\widetilde\nu$ has expression $(\ref{eq:m-range-infty-1})$, $(\ref{eq:m-h})$ and 
\begin{equation}
\label{eq:m-mean-zero}
\int_{\mathbb R^d} x\widetilde \mu(dx)=0. 
\end{equation}
\end{thm}

As seen in the proof of Lemma \ref{lem:m-eq:relation-widetilde-nu}, the function 
$\widetilde \ell_\xi(x)$ is given by 
$$
\widetilde \ell_\xi(x) = \int_x^\infty r^\alpha \nu_\xi(dr),
$$
where $\nu_\xi$ is the radial component of the L\'evy measure $\nu$ of $\mu \in
\mathfrak D(\Phi_{\al}^{m+1})$.

A function $f(t)$ defined for $t>0$ is called $m$-times monotone where $m$ is an integer, $m\ge 2$, if
$(-1)^kf^{(k)}(t)$ is nonnegative, nonincreasing and convex for $t>0$,
and for $k=0,1,2,...,m-2$.
When $m=1$, $f(t)$ will simply be nonnegative and nonincreasing.

Note that $\widetilde h_\xi^{(m)}(u)$ is  $m$-times monotone. 
In order to see this, we have only to differentiate it in the following way.
\begin{align*}
\frac{d}{ds}\widetilde h_\xi^{(m)}(s) &=
-\frac{1}{s}\int_s^\infty \frac{1}{(m-2)!x}\left(\log (x/s)
\right)^{m-2}dx \int_x^\infty r^\alpha \nu_\xi(dr)<0, \\
\frac{d^2}{ds^2}\widetilde h_\xi^{(m)} (s) &= 
\frac{1}{s^2}\int_s^\infty \frac{1}{(m-2)!x}\left(\log  (x/s) \right)^{m-2}dx
\int_x^\infty r^\alpha \nu_\xi(dr)\\
&\hskip 10mm +\frac{1}{s^2}\int_{s}^{\infty}\frac{1}{(m-3)!x}\left(\log
 (x/s) \right)^{m-3}dx
\int_{x}^{\infty}r^\alpha \nu_\xi(dr) >0.
\end{align*}
The differentiation continues to $m-1$ times, but $(d/ds)^{m-1}\widetilde h_\xi^{(m)}(s)$ includes the term
$$
(-s)^{1-m}\int_{s}^{\infty}{x}^{-1}dx\int_{x}^{\infty}r^\alpha\nu_\xi(dr) 
$$
and hence $(d/ds)^m \widetilde h_\xi^{(m)}(s)$ includes the term
$$
(-s)^{-m}\int_{s}^{\infty}r^\alpha\nu_\xi(dr).
$$
Then since we have no information about absolute continuity of the
measure $\nu_\xi(dr)$ and differentiability of $\int_s^\infty
r^\alpha\nu_\xi(dr)$ can not be guaranteed, we can not assert any
stronger results for $\widetilde h_\xi^{(m)}(s)$ other than $m$-times
differentiability.

We need the following Lemma and here we use the same notations as
before. This lemma follows from similar arguments as those used in Lemma
4.4 in Sato (2006b). 

\begin{lem}
\label{lem:m-eq:relation-widetilde-nu}
Let $-\infty<\alpha<2$ and $m=1,2,\ldots$, and let $\wt\nu$ be a L\'evy measure.
Then there exists a L\'evy measure $\nu$ satisfying $(6.8)$ such that
\begin{eqnarray}
\label{eq:m-relation-widetilde-nu}
\begin{cases}
\int_{\mathbb R^d}(|x|^2 \land 1)\nu(dx)<\infty, & \text{when}\ \alpha<0,\\
\int_{\mathbb R^d}(|x|^2 \land 1)(\log^\ast |x|)^{m+1}\nu(dx)<\infty, & \text{when}\ \alpha=0,\\
\int_{\mathbb R^d}(|x|^2 \land |x|^\alpha)(\log^\ast |x|)^m
\nu(dx)<\infty, & \text{when}\ 0<\alpha<2\\
\end{cases}
\end{eqnarray}
if and only if $\widetilde \nu$ is represented as $(\ref{eq:m-range-infty-1})$.
\end{lem}

\n
{\it Proof of Lemma \ref{lem:m-eq:relation-widetilde-nu}. }
We first prove the only if part.
Assume that the L\'evy measure $\nu$ satisfy (\ref{eq:m-nu-mapped}) and (\ref{eq:m-relation-widetilde-nu}). 
The polar decomposition gives
\begin{eqnarray}
\label{m-cond-polar-decomp}
\begin{cases}
\int_{S}\lambda(d\xi)\int_0^\infty(r^2 \land 1)\nu_{\xi}(dr)<\infty,\
 \text{when}\ \alpha< 0, \\
\int_{S}\lambda(d\xi)\int_0^\infty(r^2 \land 1)(\log^\ast r)^{m+1} \nu_{\xi}(dr)<\infty,\
 \text{when}\ \alpha= 0, \\
\int_{S}\lambda(d\xi)\int_0^\infty(r^2 \land r^\alpha)(\log^\ast r)^m \nu_{\xi}(dr)<\infty,\
 \text{when}\ \alpha> 0. 
\end{cases}
\end{eqnarray}
Then we have for $B\in\mathcal B(\mathbb R^d)$
{\allowdisplaybreaks
\begin{align*}
\widetilde \nu(B) &= (m!)^{-1}\int_0^1 t^{-\alpha-1}\nu(t^{-1}B)(\log t^{-1})^m dt\\
&=(m!)^{-1}\int_0^1 t^{-\alpha-1}dt \int_S\lambda(d\xi) \int_0^\infty 1_B(t\xi r)(\log t^{-1})^m \nu_\xi(dr) \\
&= (m!)^{-1}\int_S \lambda(d\xi) \int_0^\infty  r^\alpha\nu_\xi(dr)
\int_0^r 1_B(s\xi)s^{-\alpha-1}  (\log (r/s))^m ds \\
&=:  \int_S \lambda(d\xi) \int_0^\infty 1_B(s\xi)s^{-\alpha-1}\widetilde h^m_\xi(s)ds,   
\end{align*}
where 
\begin{align*}
\widetilde h_\xi^{(m)}(s) &= (m!)^{-1} \int_s^\infty (\log (r/s))^m r^\alpha
 \nu_\xi(dr) \\
&=  ({(m-1)!})^{-1} \int_s^\infty r^\alpha \nu_\xi(dr) \int_s^r x^{-1}(\log (x/s))^{m-1}dx \\
&=  ({(m-1)!})^{-1}\int_s^\infty {x}^{-1}(\log (x/s))^{m-1}dx \int_x^\infty r^\alpha \nu_\xi(dr) \\
&=: ({(m-1)!}) ^{-1}\int_s^\infty {x}^{-1}(\log (x/s))^{m-1}
 \widetilde \ell_\xi (x)dx. 
\end{align*}
Here $\widetilde \ell_\xi(u)$ is measurable in $\xi$ and for
$\lambda-a.e. \xi.$ nonincreasing in $u\in (0,\infty)$, and 
$\lim_{u\to\infty}\widetilde \ell_\xi(u)=0$ from
(\ref{m-cond-polar-decomp}).

Conversely, suppose that $\widetilde \nu$ satisfies
(\ref{eq:m-nu-mapped}). We consider the case $-\infty<\alpha<0$.
Then since $h_\xi^{(m)}(r)$ is a continuous decreasing function, 
we can define a Lebesgue-Stieltjes measure 
$\widetilde R_\xi$ on $(0,\infty)$ satisfying
$\widetilde R_\xi((r,s])=-\widetilde h^{(m)}_\xi(s) +\widetilde h^{(m)}_\xi(r)$ and put
$\nu_\xi(dr)=r^{-\alpha}\widetilde R_\xi(dr).$ 
Furthermore define
$$
\nu(B)=\int_S \widetilde \lambda(dr)\int_0^\infty 1_B(r\xi)\nu_\xi(dr).
$$
Here the same logic as in the proof of Lemma
\ref{lem:eq:relation-widetilde-nu} holds and
we see that (\ref{eq:m-relation-widetilde-nu}).

In the following, similar to the proof of Lemma \ref{lem:eq:relation-widetilde-nu}, we put 
$\widetilde R_\xi([r,\infty))=\widetilde \ell_\xi(r+)$
and $\nu_\xi(dr)=r^{-\alpha}\widetilde R_\xi(dr)$.
Furthermore define
$$
\nu(B) = \int_S \widetilde \lambda(d \xi)\int_0^\infty 1_B(r\xi)\nu_\xi(dr).
$$
Then for the case $\alpha=0$, let $\lambda=\widetilde \lambda$, and we have
\begin{align*}
\int_{\mathbb R^d}(|x|^2\land & 1)(\log^\ast |x|)^{m+1}\nu(dx) \\
&= \int_S\lambda(d\xi) \int_0^\infty(r^2\land 1)(\log^\ast
 r)^{m+1}\nu_\xi(dr) \\ 
& =\int_S\lambda(d\xi)\left(\int_0^1 r^{2}\widetilde R_\xi(dr)+\int_1^\infty (\log^\ast
r)^{m+1}\widetilde R_\xi(dr)\right).
\end{align*}
Since $\widetilde \nu$ is a L\'evy measure, it follows that
$$
 \int_S\lambda(d\xi)((m-1)!)^{-1}\int_0^\infty(r^2\land 1)r^{-1}dr\int_r^\infty x^{-1}
(\log (x/r))^{m-1}\widetilde \ell_\xi(x)dx<\infty.
$$ 
Then a simple calculation gives
\begin{align*}
0< ((m-1)!)^{-1}&\int_0^1 rdr \int_r^\infty x^{-1}(\log (x/r))^{m-1}dx\int_x^\infty \widetilde R_\xi(dy) \\
&=  ((m-1)!)^{-1}\int_0^1 rdr \int_r^\infty \widetilde R_\xi(dy)\int_r^y
x^{-1}(\log (x/r))^{m-1}dx \\
& = (m!)^{-1}\int_0^1 rdr \int_r^\infty (\log (y/r))^m \widetilde R_\xi(dy) \\
& = (m!)^{-1}\int_0^1  \widetilde R_\xi(dy)  \int_0^y r(\log (y/r))^m  dr \\
& \hskip 10mm +(m!)^{-1}\int_1^\infty \widetilde R_\xi(dy)  \int_0^1
 r(\log (y/r))^m  dr<\infty.  
\end{align*}
Since the last two integrals are positive and 
\begin{align*}
(m!)^{-1}&\int_0^1 \widetilde R_\xi(dy) \int_0^y  r(\log (y/r))^m dr  \\
& =(m!)^{-1}\int_0^1 t(\log t^{-1})^m dt \int_0^1  y ^2 \widetilde R_\xi(dy),
\end{align*}
the finiteness of $\int_S\lambda(d\xi)\int_0^1r^2 \widetilde R_\xi(dr)$ is shown. 
Next we see that
\begin{align*}
0<((m-1)!)^{-1}&\int_1^\infty r^{-1}dr  \int_r^\infty x^{-1}(\log (x/r))^{m-1}dx
\int_x^\infty \widetilde R_\xi(dy) \\
&=((m-1)!)^{-1}\int_{1}^{\infty}r^{-1}dr \int_r^\infty \widetilde R_\xi(dy)\int_r^y x^{-1}(\log (x/r))^{m-1}dx \\
&=(m!)^{-1}\int_1^\infty r^{-1}dr \int_r^\infty \left[ (\log (x/r))^m \right]_r^y \widetilde R_\xi(dy) \\
&=(m!)^{-1}\int_1^\infty \widetilde R_\xi(dy) \int_1^y r^{-1}(\log (y/r))^m dr \\
& =((m+1)!)^{-1}\int_1^\infty (\log y)^{m+1}\widetilde R_\xi(dy)<\infty.
\end{align*}
Hence, we have (\ref{eq:m-relation-widetilde-nu}).

When $0<\alpha<2$, 
let $\lambda=\widetilde\lambda$ and we have
\begin{align*}
\int_0^\infty & (|x|^2\land |x|^\alpha)(\log^\ast |x|)^m \nu(dx) \\
& =\int_S\lambda(d\xi)\int_0^\infty(r^2\land r^\alpha)(\log^\ast r)^m\nu_\xi(dr)\\
& =\int_S\lambda(d\xi)\int_0^\infty (r^2\land r^\alpha)(\log^\ast r)^mr^{-\alpha}\widetilde R_\xi(dr)\\
& =\int_S\lambda(d\xi)\left(\int_0^1 r^{2-\alpha}\widetilde R_\xi(dr)+\int_1^\infty (\log^\ast r)^m
\widetilde R_\xi(dr)\right).
\end{align*}
Since $\widetilde \nu$ is a L\'evy measure. We have
$$
\int_S\lambda(d\xi)((m-1)!)^{-1}\int_0^\infty(r^2\land 1)r^{-\alpha-1} dr \int_r^\infty
x^{-1}(\log (x/r))^{m-1}\widetilde \ell_\xi(x)dx<\infty.
$$ 
Thus 
\begin{align*}
0<((m-1)!)^{-1}& \int_0^1 r^{1-\alpha}dr \int_r^\infty x^{-1}(\log (x/r))^{m-1}dx\int_x^\infty  
\widetilde R_\xi(dy)\\
& =((m-1)!)^{-1}\int_0^1r^{1-\alpha}dr \int_r^\infty \widetilde R_\xi(dy)\int_r^y x^{-1}(\log (x/r))^{m-1}dx\\
& =(m!)^{-1} \int_0^1  \widetilde R_\xi(dy)  \int_0^y r^{1-\alpha}(\log
 (y/r))^m dr \\
& \hskip 10mm +(m!)^{-1}\int_1^\infty \widetilde R_\xi(dy) \int_0^1
 r^{1-\alpha}(\log (y/r))^m dr<\infty. 
\end{align*}
The first term in the right-hand side equals
$$
(m!)^{-1}\int_0^1 y^{2-\alpha} \widetilde R_\xi(dy)\int_0^1 t^{1-\alpha}(\log
t^{-1})^m dt.
$$
Furthermore,
\begin{align*}
0<((m-1)!)^{-1}&\int_1^\infty r^{-\alpha-1}dr\int_r^\infty x^{-1}(\log (x/r))^{m-1}\widetilde \ell_\xi(x)dx \\
& =((m-1)!)^{-1}\int_1^\infty r^{-\alpha-1}dr \int_r^\infty x^{-1} (\log (x/r))^{m-1}dx
\int_x^\infty \widetilde R_\xi(dy) \\
& =((m-1)!)^{-1}\int_1^\infty r^{-\alpha-1}dr \int_r^\infty \widetilde R_\xi(dy)\int_r^y x^{-1}(\log (x/r))^{m-1}dx \\
& =(m!)^{-1}\int_1^\infty r^{-\alpha-1}dr   \int_r^\infty (\log (y/r))^m \widetilde R_\xi(dy)\\
& =(m!)^{-1}\int_1^\infty \widetilde R_\xi(dy) \int_1^y r^{-\alpha-1}(\log (y/r))^m dr<\infty. 
\end{align*}
Here with the integral by parts formula 
$$
 (m!)^{-1} \int_1^y r^{-\alpha-1} (\log (y/r))^mdr
$$
is a linear combination of 
$(\log y)^k,\ k=0,\ldots,m$, and the coefficient of $(\log y)^m$ is positive.
Thus, (\ref{eq:m-relation-widetilde-nu}) holds.
\hfill $\Box$

\vskip 3mm
\noindent
{\it Proof of Theorem \ref{thm:range-of-K-m-mapping}.}\\
$(i)\ (-\infty<\alpha<1)$ \\
Assume that $\widetilde \mu\in\Phi_{\al}^{m+1}(\frak D(\Phi_{\al}^{m+1}))$. The Gaussian
case is obvious. When $\widetilde \mu$ is non-Gaussian, then from Lemma
\ref{lem:m-eq:relation-widetilde-nu} there exists $\widetilde \nu$ satisfying (\ref{eq:m-range-infty-1})
and (\ref{eq:m-h}). 

We see the converse. If $\widetilde \mu$ is Gaussian,
then putting $A=(2-\alpha)^m \widetilde A,\ \nu=0$ and
$\gamma=(1-\alpha)^m\widetilde \gamma$, we have
$\mu=\mu_{(A,\nu,\gamma)}\in \mathfrak D(\Phi_{\al}^{m+1})$ and $\widetilde
\mu=\Phi_{\al}^{m+1}(\mu)$. 

If $\widetilde \mu$ is non-Gaussian, then (\ref{eq:m-range-infty-1})
and (\ref{eq:m-h}) give the measure $\nu$ in Lemma
\ref{lem:m-eq:relation-widetilde-nu}. We put $A=(2-\alpha)^m\widetilde A$ and 
$$
\gamma=(1-\alpha)^m\left(
\widetilde \gamma +(m!)^{-1}\int_0^1 s^{-\alpha}(\log s^{-1})^m ds \int_{\mathbb
R^d}x\left(
\frac{1}{1+|x|^2}-\frac{1}{1+|sx|^2}
\right)\nu(dx)
\right)
$$ 
The existence of $\gamma$ is proved as follows. Let $c_1$ and $c_2$ be some
positive constants. Then
\begin{align*}
\int_0^1 s^{-\alpha}&(\log s^{-1})^mds \int_{\mathbb R^d}|x|\left|
\frac{1}{1+|x|^2}-\frac{1}{1+s^2|x|^2}\right|\nu(dx)\\
&\le \int_0^1 s^{-\alpha}(\log s^{-1})^m ds \left(\int_{|x|\le 1}\frac{|x|^3}
{(1+|x|^2)(1+s^2|x|^2)}\nu(dx) \right. \\
& \left. \hskip 10mm +c_1\int_{|x|>1,|sx|\le 1}|x|\nu(dx)+c_2\int_{|x|>1,|sx|>1} s^{-2}|x|^{-1}\nu(dx)  \right) \\
&= \int_0^1 s^{-\alpha}(\log s^{-1})^m ds \int_{|x|\le 1}|x|^3 \nu(dx) \\
& \quad \hskip 10mm +c_1\int_{|x|>1}|x|\nu(dx)\int_0^{1/|x|}s^{-\alpha}(\log s^{-1})^m ds \\
& \quad \hskip 10mm+c_2\int_{|x|>1}|x|^{-1}\nu(dx)\int_{1/|x|}^1 s^{-\alpha-2}(\log s^{-1})^m ds.
\end{align*}
Here by the  integral by parts formula, $\int_0^{1/|x|}s^{-\alpha}(\log s^{-1})^m ds$ is shown to be constructed 
by a linear combination of $|x|^{\alpha-1}(\log |x|)^k$ with $k=0,1,\ldots,m$ and 
$\linebreak$ 
$\int_{1/|x|}^{1}s^{-\alpha-2}(\log s^{-1})^m ds$ is shown to be constructed by a linear combination of 
$\linebreak$ 
$|x|^{\alpha+1}(\log |x|)^k$ with $k=0,1,\ldots,m$. 
Then on behalf of (\ref{eq:m-relation-widetilde-nu}) we can prove the existence of
$\gamma$. Thus $\mu=\mu_{(A,\nu,\gamma)}\in \mathfrak D(\Phi_{\al}^{m+1})$
and $\Phi_{\al}^{m+1}(\mu)=\widetilde \mu$. \\
$(ii)\ (\alpha=1)$ \\
Suppose that $\widetilde \mu\in \Phi_1^{m+1}(\mathfrak D(\Phi_{1}^{m+1}))$ and
$\mu=\mu_{(A,\nu,\gamma)}\in\mathfrak D(\Phi_1^{m+1})$. First assume that
$\widetilde \mu$ is Gaussian. Then for given $\varphi\in C_{\#}^+$ (see
the beginning of Proof of Theorem 4.7 for its definition),
$$
0=\int_0^1 s^{-2}(\log s^{-1})^m ds \int_{\mathbb R^d}\varphi(s^{-1}x)\nu(dx),
$$
which implies $0=s^{-2}(\log s^{-1})^m \int_{\mathbb
R^d}\varphi(s^{-1}x)\nu(dx).$ Since by the dominated convergence theorem 
$s^{-2}(\log s^{-1})^m \int_{\mathbb R^d}\varphi(s^{-1}x)\nu(dx)$
is continuous in $s$, letting $s=1/2$, we have $\nu=0$. 
This together with
$\gamma=0$ (which follows from Proposition 4.3) implies $\widetilde
\gamma=0$. Hence $\widetilde \mu$ is centered Gaussian.
If $\widetilde \mu$ is non-Gaussian, then Lemma
\ref{lem:m-eq:relation-widetilde-nu} assures the existence of a measure $\widetilde \nu$
such that satisfies (\ref{eq:m-range-infty-1})
and (\ref{eq:m-h}) with $\alpha=1$, and the combination of 
$\gamma=-\int_{\mathbb R^d}\frac{x|x|^2}{1+|x|^2}\nu(dx)$ and (\ref{eq:m-gamma-mapped}) implies
(\ref{eq:cond-conve-gamma-m-mapped}). 

We prove the converse part. If $\widetilde \mu$ is centered Gaussian,
then $\widetilde \mu\in \Phi_{\al}^{m+1}(\mathfrak D(\Phi_{\al}^{m+1}))$ by Theorem \ref{range-m0}. If
$\widetilde \mu$ is non-Gaussian and satisfies (\ref{eq:m-range-infty-1})
and (\ref{eq:m-h}), then  Lemma \ref{lem:m-eq:relation-widetilde-nu} yields
the existence of a measure $\nu$ satisfying (\ref{eq:m-nu-mapped}) and 
(\ref{eq:m-relation-widetilde-nu}) with $\alpha=1$. Let $\gamma=-\int_{\mathbb
R^d}\frac{x|x|^2}{1+|x|^2}\nu(dx), A=\widetilde A$ and
$\mu=\mu_{(A,\nu,\gamma)}$. It follows from
(\ref{eq:cond-conve-gamma-m-mapped}) that
$$
\lim_{T\to\infty}\int_0^T \varepsilon_{\alpha,m}^\ast(t)dt \int_{\mathbb
R^d} \frac{x|\varepsilon_{\alpha,m}^\ast(t)
x|^2}{1+|\varepsilon_{\alpha,m}^\ast(t)x|^2}\nu(dx)\quad \text{exists
in}\ \mathbb{R}^d.
$$
Then Proposition 2.6 (ii) of Sato (2006a) assures $\mu\in \mathfrak D(\Phi_{\al}^{m+1})$. 
Furthermore (\ref{eq:cond-conve-gamma-m-mapped}) implies 
$$
\widetilde \gamma =\lim_{T\to\infty}\int_0^T
\varepsilon_{\alpha,m}^\ast(t)dt 
\left(
\gamma-\int_{\mathbb R^d}x\left(
\frac{1}{1+|x|^2}-\frac{1}{1+|\varepsilon_{\alpha,m}^\ast(t)x|^2}
\right)\nu(dx)
\right).
$$
Thus looking at (\ref{eq:m-gamma-mapped}), we confirm that $\Phi_1^{m+1}(\mu)=\widetilde \mu$ and
$\widetilde \mu\in \Phi_1^{m+1}(\mathfrak D (\Phi_1^{m+1}))$. \\
(iii) $(0<\alpha<2)$ \\
Assume that $\widetilde \mu=\Phi_{\al}^{m+1}(\mu)$ with some
$\mu=\mu_{(A,\nu,\gamma)}\in \mathfrak D(\Phi_{\al}^{m+1})$. 
The Gaussian case is the same as that of Proof for (ii). 
If $\widetilde \mu$ is
non-Gaussian, then it follows from Lemma
\ref{lem:m-eq:relation-widetilde-nu} that there exists $\widetilde \nu$
satisfying (\ref{eq:m-range-infty-1}) and (\ref{eq:m-h}). Since $\mu\in
\mathfrak D(\Phi_{\al}^{m+1}),\ \nu$ and $\gamma$ satisfy
$\int_{|x|>1}|x|^\alpha(\log |x|)^m \nu(dx)<\infty$ and
$\gamma=-\int_{\mathbb R^d}\frac{x|x|^2}{1+|x|^2}\nu(dx)$, respectively. 
We show the existence of $\widetilde \gamma$, we have 
\begin{align*}
({m!})^{-1}\int_0^1 t^{2-\alpha}&(\log t^{-1})^m dt\int_{\mathbb R^d}\frac{|x|^3}{1+t^2|x|^2}\nu(dx) \\
& \le ({m!})^{-1}\int_0^1 t^{2-\alpha}(\log t^{-1})^mdt\int_{|x|\le 1}|x|^3\nu(dx)\\
&\hskip 10mm + ({m!})^{-1}\int_{|x|>1}|x|^3\nu(dx)\int_0^{1/|x|}t^{2-\alpha}(\log t^{-1})^m dt\\
&\hskip 10mm+ ({m!})^{-1}\int_{|x|>1}|x|\nu(dx)\int_{1/|x|}^1 t^{-\alpha}(\log t^{-1})^m dt.
\end{align*}
Here by the  integral by parts formula, $\int_0^{1/|x|}t^{2-\alpha}(\log
t^{-1})^m dt$ is shown to be a linear combination of 
$|x|^{\alpha-3}(\log|x|)^k$
with $k=0,1,\ldots,m$ and $\int_{1/|x|}^1t^{-\alpha}(\log t^{-1})^m dt$
is shown to be a linear combination of
$|x|^{\alpha-1}(\log |x|)^k$
with $k=0,1,\ldots,m$. Then form (\ref{eq:m-relation-widetilde-nu})
$\widetilde \gamma$ exists and equals to 
\begin{align*}
\widetilde \gamma &= -\int_0^\infty
 \varepsilon_{\alpha,m}^\ast(t)dt\int_{\mathbb
 R^d} \frac{x|\varepsilon_{\alpha,m}^\ast(t) x|^2 }
 {1+|\varepsilon_{\alpha,m}^\ast(t) x|^2}\nu(dx)\\
&= -\int_{\mathbb R^d}\frac{x|x|^2}{1+|x|^2}\widetilde \nu(dx),
\end{align*} 
which is (\ref{eq:m-mean-zero}).

We consider the converse. 
If $\widetilde \mu$ is centered Gaussian with its component $\widetilde
A$. Then by putting $A=(2-\alpha)^{m+1}\widetilde A$, we have
$\mu_{(A,0,0)}\in \mathfrak D(\Phi_\alpha^{m+1})$ and
$\widetilde \mu = \Phi_\alpha^{m+1}(\mu)$. 
Suppose
$\widetilde \mu$ be non-Gaussian and satisfy condition of (iii). 
On behalf of Lemma \ref{lem:m-eq:relation-widetilde-nu}, we have a measure
$\nu$ satisfying (\ref{eq:m-nu-mapped}) and (\ref{eq:m-relation-widetilde-nu}). 
We investigate the absolute moment of $\widetilde \nu$ and see that
\begin{align*}
\int_{\mathbb R^d}\frac{|x|^3}{1+|x|^2}\widetilde\nu(dx) &=
\int_0^\infty \varepsilon_{\alpha,m}^\ast(t)dt \int_{\mathbb R^d}
\frac{|x|^3}{1+|\varepsilon_{\alpha,m}^\ast(t)x|^2}\nu(dx) \\
&=\frac{1}{m!}\int_0^1 s^{2-\alpha}(\log s^{-1})^m ds \int_{\mathbb R^d}
 \frac{|x|^3}{1+s^2|x|^2}\nu(dx). 
\end{align*}
Then as we have seen in the preceding paragraph,
$$
\int_{\mathbb R^d}\frac{|x|^3}{1+|x|^2}\widetilde \nu(dx)<\infty.$$
Thus $\int_{\mathbb
R^d}|x|\widetilde \mu(dx)<\infty$, and hence
$\int_{\mathbb R^d} x \widetilde \mu(dx)=0$. 
Let $\gamma=-\int_{\mathbb{R}^d}\frac{x|x|^2}{1+|x|^2}\nu(dx)$ and
$A=(2-\alpha)^{m+1}\widetilde A$. Then on behalf of
(\ref{eq:m-nu-mapped}), Theorem
\ref{thm:integral-representation-domain-K} $(\mathrm{v})$ is satisfied. Thus
$\mu=\mu_{(A,\nu,\gamma)}\in \mathfrak D(\Phi_{\alpha}^{m+1})$. 
\hfill $\Box$
}


\vskip 10mm
\addcontentsline{toc}{section}{References}

\begin{thebibliography}{99}


\bibitem{1}  O.E. Barndorff-Nielsen and M. Maejima  (2008).
Semigroups of Upsilon transformations, to appear in {\it Stoch. Proc. Appl.}

\bibitem{1}  O.E. Barndorff-Nielsen, M. Maejima and K. Sato (2006).
Some classes of multivariate infinitely divisible distributions 
admitting stochastic integral representations, {\it Bernoulli},
{\bf 12}, 1--33.

\bibitem{1}
Z.J. Jurek (1983). The class $L_m(Q)$ of probability measures on Banach spaces,
{\it Bull. Polish Acad. Sci. Math.} {\bf 31}, 51--62.


\bibitem{3} Z.J. Jurek (1985). Relations between the 
{\it s}--selfdecomposable and selfdecomposable measures, {\it Ann. Probab.} {\bf 13}, 
592--608.

\bibitem{1} Z.J. Jurek (1988).
Random integral representations for classes of limit distributions similar to Levy class $L_0$.
{\it Probab. Th. Rel. Fields} {\bf 78}, 473--490.

\bibitem{1} Z.J. Jurek (2004).
The random integral representation hypothesis revisited : new class of s-selfdecomposable laws,
in : {\it Abstract and Applied Analysis} ; Proc. Intern. Conf. Hanoi, 2002, World Scientific, 495--514.

\bibitem{1} 
T.A. O'Connor (1979). Infinitely divisible distributions similar to class L distributions, {\it Z.W.} {\bf 50}, 265--271.)

\bibitem{1}
J. Rosi\'nski (2007).  Temparing stable processes,
{\it Stoch. Proc. Appl.} {\bf 117}, 677--707.
\bibitem{1}
K. Sato (1999).  {\it L\'evy Processes and Infinitely Divisible Distributions}. 
Cambridge Univ.\ Press, Cambridge.

\bibitem{4} K. Sato (2006a). Additive processes and stochastic integrals,
{\it Illinois J. Math.} {\bf 50} (Doob Volume), 825--851. 


\bibitem{5} K. Sato (2006b). Two families of improper stochastic integrals with respect to L\'evy processes,
{\it ALEA Lat. Am. J.
Probab. Math. Statist.} {\bf 1}, 47--87.

\bibitem{1} S.J. Wolfe (1982). On a continuous analogue of the stochastic difference equation 
$X_n = \rho X_{n-1}+ B_n$, {\it Stoch. Proc. Appl.} {\bf 12}, 301--312. 

\end{thebibliography}
\end{document}